\documentclass[11pt,a4paper]{article}

\usepackage{mathtools}
\usepackage{amssymb}
\usepackage{amsthm}

\usepackage[T1]{fontenc}
\usepackage[utf8]{inputenc}
\usepackage[english]{babel}

\usepackage{graphicx}
\usepackage{xcolor}

\usepackage{tikz}
\usetikzlibrary{arrows.meta, positioning, calc}

\usepackage[margin=3cm]{geometry}
\usepackage[hidelinks]{hyperref}
\usepackage[nameinlink,capitalise,noabbrev]{cleveref}
\usepackage[protrusion=true,expansion=false]{microtype}
\usepackage{float}

\numberwithin{equation}{section}

\theoremstyle{plain}
\newtheorem{theorem}{Theorem}[section]
\newtheorem{proposition}[theorem]{Proposition}
\newtheorem{lemma}[theorem]{Lemma}
\newtheorem{corollary}[theorem]{Corollary}
\newtheorem{assumption}[theorem]{Assumption}

\theoremstyle{definition}
\newtheorem{definition}[theorem]{Definition}
\newtheorem{example}[theorem]{Example}
\newtheorem{question}[theorem]{Open question}

\theoremstyle{remark}
\newtheorem{remark}[theorem]{Remark}

\crefname{theorem}{Theorem}{Theorems}
\crefname{proposition}{Proposition}{Propositions}
\crefname{lemma}{Lemma}{Lemmas}
\crefname{corollary}{Corollary}{Corollaries}
\crefname{definition}{Definition}{Definitions}
\crefname{remark}{Remark}{Remarks}
\crefname{assumption}{Assumption}{Assumptions}
\crefname{question}{Open question}{Open questions}
\crefname{equation}{Eq.}{Eqs.}
\crefname{section}{Section}{Sections}
\crefname{subsection}{Section}{Sections}

\newcommand{\R}{\mathbb{R}}
\newcommand{\N}{\mathbb{N}}
\newcommand{\wA}{\omega_A}
\newcommand{\Astar}{A^{\ast}}

\newcommand{\Cb}{C_b}
\newcommand{\Cz}{C_0}

\newcommand{\Msp}{\mathcal{M}}
\newcommand{\Xdec}{\mathcal{X}}
\DeclareMathOperator{\supp}{supp}
\DeclareMathOperator{\osc}{osc}
\DeclareMathOperator{\Lip}{Lip}

\begin{document}

\title{Weighted Sup-Norms at Infinity:\\
Exhaustion Functions, Decay Classes, and Anisotropic Ends}
\author{Armen Petrosyan}
\date{\today}
\maketitle

\begin{abstract}
Let $A$ be a non-compact, locally compact, $\sigma$-compact Hausdorff space and let
$h:A\to[0,\infty)$ be a continuous proper exhaustion function. For an admissible weight
$\varphi$ we study the weighted sup-norms
$\|f\|_{\infty,h;\varphi;L}=\sup_A \varphi(h)\,|f-L|$ and their quotient (``sharp'')
version modulo constants, which measure \emph{how fast} a function approaches its limit
at infinity.

Our main structural result is a reduction theorem: when $\varphi$ is unbounded, the
infimum over constants defining the sharp norm is attained at the unique constant
$L=\lim_{a\to\wA}f(a)$, and the map $[f]\mapsto \varphi(h)\,(f-L)$ is a surjective linear
isometry from the sharp space onto $\Cb(A)$, restricting to an isometry of the decay class
onto $\Cz(A)$. Consequently the isotropic theory---completeness, duality, compactness,
density of truncations---is exactly the classical theory of $\Cb(A)$ and $\Cz(A)$ read
through this isometry. We spell out what this does and does not give: in particular the
natural space of weighted zero-mass measures $\Msp_0^{W^{-1}}(A)$ embeds isometrically
into the dual of the decay class with dense range but is \emph{not} the whole dual, and
norming measures generally fail to exist.

We also show by counterexample that continuous proper exhaustion functions are \emph{not}
in general coarsely affinely equivalent, so the decay classification $O(h^{-p})$ genuinely
depends on the choice of $h$ within its coarse class; all our statements are therefore made
relative to a fixed coarse-affine class, under which they are invariant.

Finally we treat anisotropic spaces with finitely many ends. We show that the naive
single-constant quotient is the wrong object there---it is infinite as soon as two ends
carry different limits---and replace it by an end-by-end functional for which a
gluing lemma, a structure theorem and completeness hold.

An appendix lists, statement by statement, the corrections made with respect to the
posted version of this work; a number of results announced there are false as stated and
are retracted here.
\end{abstract}

\tableofcontents

\section{Introduction}
\label{sec:intro}

\subsection{Motivation}

In many parts of analysis one wants to say not only \emph{that} a function converges at
infinity but \emph{how fast}. The classical device is a weighted sup-norm: one fixes a
weight $w$ growing at infinity and measures $\sup_a w(a)\,|f(a)-L|$. Finiteness of this
quantity is exactly a decay statement, $f-L=O(1/w)$.

The weight is usually chosen by hand. On $\R^n$ or on a Riemannian manifold the natural
choice is a power of the distance to a base point; on a manifold with conic or cylindrical
ends one uses the weights of the $b$-calculus \cite{Melrose1993, Kondratev1967, Agmon1965}.
The purpose of this note is to organise this practice around a single object, the
\emph{exhaustion function} $h$, and to determine precisely how much structure the resulting
scale of spaces actually carries.

The answer, in the isotropic case, is: exactly the structure of $\Cb(A)$ and $\Cz(A)$, and
nothing more. We regard this as the main point of the paper rather than as a defect: it
identifies the weighted sharp spaces up to isometry, it makes every functional-analytic
question about them immediately answerable, and it delimits sharply where genuinely new
phenomena can occur (bounded weights, and anisotropic ends).

\subsection{Setting and main results}

Throughout, $A$ is a non-compact, locally compact, $\sigma$-compact Hausdorff space and
$h:A\to[0,\infty)$ a continuous proper exhaustion function (\cref{lem:proper-exhaustion}).
We write $\wA$ for the point at infinity of the Alexandroff compactification and
$W:=\varphi\circ h$ for the weight attached to an admissible comparison function $\varphi$
(\cref{def:admissible}).

\paragraph{(1) Failure of coarse-affine equivalence.} It is tempting to believe that the
decay classification $O(h^{-p})$ is independent of the choice of exhaustion. It is not.
On $A=\R$, the functions $h(x)=|x|$ and $h'(x)=x^{2}$ are both continuous and proper, they
are not coarsely affinely equivalent, and $O(h^{-1})\ne O(h'^{-1})$
(\cref{ex:no-coarse-affine}). What is true, for a \emph{fixed} coarse-affine class, is that the norms attached to two
coarsely affinely equivalent exhaustions are equivalent for polynomial weights
(\cref{lem:equiv-affine}). For general weights we characterise invariance
(\cref{thm:coarse-stability}). Two properties must be distinguished: comparability of the
asymptotic constants for all $h'\sim h$, and mere preservation of their vanishing. The
first holds if and only if $\varphi(\lambda t)\asymp\varphi(t)$ uniformly for $t$ in the
\emph{image} $h(A)$, for every $\lambda>0$; the second, if and only if the one-sided half
of that condition holds. Both reduce to global $O$-regular variation of $\varphi$ when
$h(A)$ contains a half-line, e.g.\ when $A$ is connected, but not otherwise: on $A=\N$ we
exhibit a weight with unbounded global dilation for which invariance nevertheless holds
(\cref{ex:N-inv-strict}), and one for which vanishing is preserved but comparability fails
(\cref{ex:N-van-strict}). Polynomial and log-polynomial weights always qualify
(\cref{cor:poly-invariance}); exponential weights do not on a connected space
(\cref{ex:exp-not-invariant}). We fix a coarse-affine class once and for all
(\cref{ass:standing}).

\paragraph{(2) The reduction theorem.} Assume $\varphi$ is unbounded, so that
$W(a)\to\infty$ as $a\to\wA$. Then (\cref{thm:reduction}) for $f\in C(A)$:
\[
\inf_{c\in\R}\ \sup_{a\in A} W(a)\,|f(a)-c| <\infty
\iff
\text{$f$ has a limit $L$ at $\wA$ and $W\,(f-L)$ is bounded},
\]
the infimum is attained at $c=L$ and at no other constant, and
\[
T:\ [f]\longmapsto W\cdot\bigl(f-L(f)\bigr)
\]
is a \emph{surjective linear isometry} onto $\Cb(A)$. It maps the decay class
\[
\Xdec:=\bigl\{[f]:\ \textstyle\limsup_{h\to\infty}W|f-L|=0\bigr\}
\]
isometrically onto $\Cz(A)$ (\cref{cor:isometry}).

Everything in \cref{sec:norms,sec:duality} is then a transcription. Completeness is the
completeness of $\Cb(A)$; density of compactly supported truncations is density of
$C_c(A)$ in $\Cz(A)$; the weighted Arzel\`a--Ascoli theorem is the classical compactness
criterion in $\Cz(A)$; the dual of $\Xdec$ is $\Msp(A)$ by Riesz representation.

\paragraph{(3) Duality is not what one expects.} The natural candidate for the dual of
$\Xdec$ is the space $\Msp_0^{W^{-1}}(A)$ of finite signed Radon measures of \emph{zero
total mass} with $\int W^{-1}d|\nu|<\infty$, acting by $[f]\mapsto\int f\,d\nu$. We show
(\cref{thm:duality}) that $\nu\mapsto\Lambda_\nu$ is an isometric embedding with
\emph{norm-dense range}, but that it is \emph{not surjective}: its image is
$\{\sigma\in\Msp(A):\int W\,d|\sigma|<\infty,\ \int W\,d\sigma=0\}$ under the Riesz
identification, and $\delta_{a_0}$ is not in it. As a consequence, norming measures need
not exist: we exhibit $[f]$ for which the supremum over the unit ball of
$\Msp_0^{W^{-1}}(A)$ is not attained (\cref{ex:no-extremal}).

\paragraph{(4) Anisotropic ends.} If $A=K_0\sqcup U_1\sqcup\dots\sqcup U_m$ with a
compact core and finitely many ends, each carrying its own exhaustion $h_i$ and weight
$\varphi_i$, the single-constant quotient is the wrong object: as soon as two ends carry
different limits and the weights are unbounded, $\inf_c\sup_A W|f-c|=+\infty$. We replace
it by the end-by-end functional
\[
\|[f]\|_{\sharp}:=\max\Bigl\{\inf_{c}\sup_{K_0}|f-c|,\ \max_i \inf_{c}\sup_{U_i}\varphi_i(h_i)|f-c|\Bigr\},
\]
which is a genuine norm on the appropriate quotient, and for which we prove a structure
theorem and completeness (\cref{sec:ends}).

\subsection{What this paper does not claim}

We make no claim of novelty for the use of weighted norms, for the construction of $h$
(this is standard point-set topology, \cite[Ch.~4]{Lee2011}), or for the recovery of the
Alexandroff compactification. The construction of a canonical weight from an exhaustion
does not, by itself, yield a decay classification independent of choices; see (1) above.
Our contribution is the identification (2), the corrected duality statement (3), and the
corrected anisotropic formalism (4), together with an explicit inventory of what fails.

\subsection{Relation to the posted version, and to the companion note}

The posted version of this work announced a number of results that are false as stated:
coarse-affine equivalence of arbitrary proper exhaustions, log-convexity of the sharp norm
in $p$, a real-interpolation theorem depending on it, a Moreau-envelope regularisation
statement, generic uniqueness of sharp minimisers, full duality with
$\Msp_0^{W^{-1}}(A)$, weak-$*$ compactness of its unit ball, existence of dual extremals,
a Banach-lattice structure, and a gluing lemma in the anisotropic setting. These are
retracted; \cref{app:corrections} lists them individually together with the reason.

The earlier version also announced, in a companion note, that the weighted Kolmogorov
metric built on this framework recovers the rate $n^{-1/2}$ in the Berry--Esseen problem
under $\mathbb{E}|X|^{2+\delta}<\infty$. That claim is withdrawn. The admissible region
for the parameters of that construction is empty unless the third moment essentially
exists; moreover the Berry--Esseen defect under heavy tails is localised at bounded
argument, where a weight that decays in $|t|$ is inactive. The companion note is being
rewritten as a negative result.

\paragraph{Notation.} $C(A)$ is the space of real continuous functions on $A$;
$\Cb(A)$ the bounded ones with the sup-norm; $\Cz(A)$ those vanishing at infinity;
$C_c(A)$ those with compact support. $\Msp(A)$ denotes the finite signed regular Borel
measures on $A$ with the total variation norm. For $S\subset A$ and $f:A\to\R$ we write
$\osc_S(f):=\inf_{c\in\R}\sup_{S}|f-c|$.

\begin{table}[ht]
\centering
\caption{Main symbols.}
\begin{tabular}{l l}
\hline
Symbol & Meaning \\
\hline
$A$ & non-compact, LCH, $\sigma$-compact space \\
$h:A\to[0,\infty)$ & continuous proper exhaustion (\cref{lem:proper-exhaustion}) \\
$\wA$, $\Astar$ & point at infinity, one-point compactification \\
$\varphi\in\Phi_{\mathrm{adm}}$ & admissible comparison function (\cref{def:admissible}) \\
$W=\varphi\circ h$ & the weight \\
$\|f\|_{\infty,h;\varphi;L}$ & fixed-$L$ weighted sup-norm \\
$\|[f]\|^{\sharp}_{\infty,h;\varphi}$ & sharp quotient norm modulo constants \\
$C^{(h)}_{\varphi}(f;L)$ & asymptotic constant \\
$\Xdec$ & decay class $\{[f]:C^{(h)}_{\varphi}([f])=0\}$ \\
$T$ & the isometry $[f]\mapsto W(f-L(f))$ (\cref{cor:isometry}) \\
$K_0,U_i,h_i,\varphi_i$ & core, ends, per-end exhaustions and weights (\cref{sec:ends}) \\
\hline
\end{tabular}
\end{table}

\section{Exhaustions and the exhaustion function}
\label{sec:exhaustions}

\subsection{Standing setting}

Throughout the paper $A$ is a \emph{non-compact}, locally compact, $\sigma$-compact
Hausdorff space. Local compactness and $\sigma$-compactness give an increasing exhaustion
by compacta $(K_n)_{n\ge0}$ with $K_n\subset\operatorname{int}(K_{n+1})$ and
$A=\bigcup_n K_n$.

\begin{lemma}[Continuous proper exhaustion function]\label{lem:proper-exhaustion}
There exists a continuous $h:A\to[0,\infty)$ such that
\begin{enumerate}
\item[(i)] $h$ is proper: $h^{-1}([0,R])$ is compact for every $R<\infty$;
\item[(ii)] $h(x)\to\infty$ exactly when $x$ leaves every compact subset of $A$.
\end{enumerate}
\end{lemma}

\begin{proof}
Set $K_{-1}:=\varnothing$. By Urysohn's lemma on a locally compact Hausdorff space there
are continuous $\phi_n:A\to[0,1]$ with $\phi_n\equiv0$ on $K_n$ and $\phi_n\equiv1$ on
$A\setminus\operatorname{int}(K_{n+1})$. Put $h:=\sum_{n\ge0} n\,\phi_n$. Each $x$ lies in
some $K_N$, so $\phi_n(x)=0$ for $n\ge N$ and the sum is locally finite, hence $h$ is
well defined and continuous. If $x\notin K_{N+1}$ then $\phi_n(x)=1$ for $0\le n\le N$,
so $h(x)\ge N(N+1)/2$. Given $R<\infty$, choose $N$ with $N(N+1)/2>R$; then
$A\setminus K_{N+1}\subset\{h>R\}$, i.e.\ $h^{-1}([0,R])\subset K_{N+1}$, which is compact.
\end{proof}

\begin{definition}[Exhaustion, exhaustion function]\label{def:exhaustion}
An \emph{exhaustion} of $A$ is a nested family $\{K_r\}_{r\ge0}$ of subsets with
$K_r\subseteq K_s$ for $r<s$ and $\bigcup_{r\ge0}K_r=A$; it is \emph{proper} if
$K_R\ne A$ for every $R<\infty$. The associated \emph{exhaustion function} is
\[
h(a):=\inf\{r\ge0:\ a\in K_r\}.
\]
\end{definition}

\begin{proposition}[Level sets and escape]\label{prop:level-sets}
For $R\ge0$ one has $\{h\le R\}=K_R$ (if the $K_r$ are closed and the infimum is attained)
and $B_R:=\{h>R\}=A\setminus K_R$. A sequence $(a_n)$ eventually leaves every $K_R$
if and only if $h(a_n)\to\infty$.
\end{proposition}

\begin{example}[The four standard instances]\label{ex:instances}
\leavevmode
\begin{itemize}
\item \emph{Metric.} $(A,d)$ proper, $K_r=\overline{B}(a_0,r)$, $h(a)=d(a_0,a)$.
\item \emph{Topological.} $K_n$ compact with $K_n\subset\operatorname{int}K_{n+1}$;
$h$ as in \cref{lem:proper-exhaustion}.
\item \emph{Directed.} $A=\N$, $K_n=\{1,\dots,n\}$, $h(n)=n$.
\item \emph{Measured.} $A=[0,\infty)$, $K_R=[0,R]$, $h(x)=x$.
\end{itemize}
In each case $h$ is the obvious function; the point of the formalism is uniformity of
notation, not the production of a new object.
\end{example}

\subsection{Coarse-affine equivalence, and its failure}

\begin{definition}[Coarse-affine equivalence]\label{def:coarse-affine}
Two functions $h,h':A\to[0,\infty)$ are \emph{coarsely affinely equivalent}, written
$h\sim h'$, if there are $a_1,a_2>0$ and $b_1,b_2\ge0$ with
\[
a_1 h-b_1\ \le\ h'\ \le\ a_2 h+b_2 \qquad\text{on }A .
\]
\end{definition}

The posted version of this work asserted that \emph{any} two continuous proper exhaustion
functions are coarsely affinely equivalent. This is false, and the failure is not exotic.

\begin{example}[Proper exhaustions need not be coarsely affinely equivalent]
\label{ex:no-coarse-affine}
Let $A=\R$, $h(x)=|x|$ and $h'(x)=x^{2}$. Both are continuous and proper. If
$h'\le a_2h+b_2$ held, then $x^{2}\le a_2|x|+b_2$ for all $x$, which is false. Moreover
the induced decay classes differ: $f(x)=(1+|x|)^{-1}$ satisfies $f=O(h^{-1})$ but not
$f=O(h'^{-1})$.
\end{example}

\begin{remark}[What survives]\label{rem:what-survives}
What is true is the following weaker, and much less useful, statement. If $h,h'$ are
continuous and proper, put $g(R):=\sup_{\{h\le R\}}h'$ and $f(R'):=\sup_{\{h'\le R'\}}h$.
Both are finite (continuous functions on compacta), non-decreasing and unbounded, and
$h'\le g\circ h$, $h\le f\circ h'$, with $R\le f(g(R))$ and $R'\le g(f(R'))$. Nothing in
this forces $f,g$ to be affinely bounded: \cref{ex:no-coarse-affine} has $g(R)=R^{2}$.
The erroneous step in the posted version was to deduce affine growth from the mutual
quasi-inverse relation.
\end{remark}

\begin{assumption}[Standing assumption]\label{ass:standing}
From now on a continuous proper exhaustion function $h$ is \emph{fixed}. All norms,
spaces and decay classes below are defined relative to $h$. For the polynomial weights
$\varphi_p$ they depend only on the coarse-affine class of $h$ (\cref{lem:equiv-affine},
\cref{cor:poly-invariance}); for general weights, invariance is equivalent to a dilation
condition on $\varphi$ along $h(A)$ (\cref{thm:coarse-stability}) and \emph{fails} for
exponential weights on a connected space (\cref{ex:exp-not-invariant}).
\end{assumption}

\begin{definition}[Equivalent exhaustions]\label{def:equiv-exh}
Two exhaustions $\{K_r\}$, $\{K'_r\}$ of $A$ are \emph{equivalent} if there are
$c_1,c_2>0$ with $K_r\subseteq K'_{c_2r}$ and $K'_r\subseteq K_{c_1r}$ for all $r\ge0$.
\end{definition}

\begin{proposition}\label{prop:equiv-functions}
If $\{K_r\}$ and $\{K'_r\}$ are equivalent exhaustions with associated functions $h,h'$,
then $c_1^{-1}h\le h'\le c_2 h$ on $A$; in particular $h\sim h'$.
\end{proposition}

\begin{proof}
Let $a\in A$ and $\varepsilon>0$. By definition of the infimum, $a\in K_{h(a)+\varepsilon}
\subseteq K'_{c_2(h(a)+\varepsilon)}$, so $h'(a)\le c_2(h(a)+\varepsilon)$; letting
$\varepsilon\downarrow0$ gives $h'(a)\le c_2h(a)$. Symmetrically
$a\in K'_{h'(a)+\varepsilon}\subseteq K_{c_1(h'(a)+\varepsilon)}$ gives $h(a)\le c_1h'(a)$.
(This avoids assuming that the infimum in \cref{def:exhaustion} is attained.)
\end{proof}

\begin{lemma}[Norms are coarse-affine invariants]\label{lem:equiv-affine}
Let $h\sim h'$ and let $\varphi(t)=(1+t)^{p}$, $p>0$. Then there are $M_1,M_2>0$ with
\[
M_1\,\|f\|_{\infty,h;\varphi;L}\ \le\ \|f\|_{\infty,h';\varphi;L}\ \le\ M_2\,\|f\|_{\infty,h;\varphi;L}
\]
for all $f$ and $L$, and likewise for the sharp norms of \cref{def:norms}.
\end{lemma}

\begin{proof}
From $h'\le a_2h+b_2$ we get $1+h'\le (1+b_2)\max(1,a_2)\,(1+h)$. From $a_1h-b_1\le h'$ we
get $h\le (h'+b_1)/a_1$, hence $1+h\le \max(1,a_1^{-1})(1+b_1)\,(1+h')$. Raise to the power
$p$, multiply by $|f-L|$ and take suprema.
\end{proof}

\section{The point at infinity}
\label{sec:infinity}

\begin{definition}[Convergence to $\wA$]\label{def:conv-infty}
Set $\Astar:=A\cup\{\wA\}$. We say $a\to\wA$ along a net or sequence if, for every $R>0$,
eventually $a\in A\setminus K_R$; equivalently, along the tail filter
$\mathcal{F}_\infty=\{A\setminus K_R\}_{R>0}$.
\end{definition}

\begin{corollary}\label{cor:conv-equiv}
$a_n\to\wA$ if and only if $h(a_n)\to\infty$.
\end{corollary}

\begin{proposition}[Alexandroff compactification]\label{prop:alexandroff}
Declare $U\subset\Astar$ open if either $U\subset A$ is open in $A$, or $\wA\in U$ and
there is $R>0$ with $(A\setminus K_R)\cup\{\wA\}\subset U$. Then $\Astar$ is a compact
Hausdorff space and the topology so defined is the Alexandroff one-point compactification
of $A$.
\end{proposition}

\begin{proof}
The neighbourhoods of $\wA$ are exactly the sets $(A\setminus K)\cup\{\wA\}$ with $K$
compact: any compact $K$ is contained in some $K_R$ by properness of $h$, and each $K_R$
is compact. This is the standard base for the Alexandroff topology, and compactness and
the Hausdorff property are the classical facts for a locally compact Hausdorff $A$.
\end{proof}

\begin{theorem}[Universal property]\label{thm:universal}
Let $Y$ be a compact Hausdorff space, $F:A\to Y$ continuous, and suppose $F$ has a limit
$y_\infty\in Y$ along $\mathcal{F}_\infty$, i.e.\ for every neighbourhood $V$ of $y_\infty$
there exists $R<\infty$ with $F(\{h>R\})\subset V$. Then there is a unique continuous
$\widetilde F:\Astar\to Y$ with $\widetilde F|_A=F$, and $\widetilde F(\wA)=y_\infty$.
\end{theorem}

\begin{proof}
Uniqueness is clear since $A$ is dense in $\Astar$ ($A$ is non-compact). Define
$\widetilde F(\wA):=y_\infty$; continuity at points of $A$ is inherited from $F$, and
continuity at $\wA$ is precisely the hypothesis, the sets $\{h>R\}\cup\{\wA\}$ forming a
neighbourhood base at $\wA$ by \cref{prop:alexandroff}.
\end{proof}

\begin{remark}
The quantifiers matter. The formulation used in the posted version (``$F(\{h>R\})\subset U$
for some neighbourhood $U$ independent of $R$'') is vacuous: taking $R=0$ and $U=Y$
satisfies it for every $F$, while the conclusion fails for, say, $A=\R$, $Y=[-1,1]$,
$F=\sin$.
\end{remark}

\begin{definition}[Limit at infinity]\label{def:limit}
For $f:A\to\R$ and $L\in\R$ we write $\lim_{a\to\wA}f(a)=L$ if for every $\varepsilon>0$
there is $R$ such that $|f(a)-L|<\varepsilon$ whenever $h(a)>R$. Limits are unique because
$A$ is non-compact.
\end{definition}

\section{Weighted sup-norms and the reduction theorem}
\label{sec:norms}

\subsection{Admissible weights and the basic functionals}

\begin{definition}[Admissible comparison functions]\label{def:admissible}
$\Phi_{\mathrm{adm}}$ is the set of $\varphi:[0,\infty)\to[1,\infty)$ such that
\begin{enumerate}
\item[(A1)] $\varphi$ is continuous, non-decreasing, $\varphi(0)=1$;
\item[(A2)] there is $K\ge1$ with $\varphi(r+s)\le K\varphi(r)\varphi(s)$ for all $r,s\ge0$.
\end{enumerate}
We call $\varphi$ \emph{unbounded} if $\varphi(t)\to\infty$ as $t\to\infty$. Standard
examples, all unbounded for the indicated parameters:
\[
\varphi_p(t)=(1+t)^{p}\ (p>0),\qquad
\varphi_{p,q}(t)=(1+t)^{p}(1+\log(1+t))^{q},\qquad
\varphi_\alpha(t)=e^{\alpha t}\ (\alpha>0).
\]
We write $\varphi_1\preceq\varphi_2$ if $\varphi_1\le C\varphi_2$ on $[R,\infty)$ for some
$C\ge1$, $R<\infty$.
\end{definition}

\begin{lemma}\label{lem:preceq-global}
If $\varphi_1\preceq\varphi_2$ then $\varphi_1\le C'\varphi_2$ on all of $[0,\infty)$,
with $C'=\max(C,\varphi_1(R))$.
\end{lemma}

\begin{proof}
On $[R,\infty)$ this is the hypothesis. On $[0,R]$, $\varphi_1\le\varphi_1(R)$ by
monotonicity while $\varphi_2\ge1$.
\end{proof}

\begin{definition}[The functionals]\label{def:norms}
Fix $\varphi\in\Phi_{\mathrm{adm}}$ and set $W:=\varphi\circ h$. For $f\in C(A)$,
$c\in\R$:
\[
\|f\|_{\infty,h;\varphi;c}:=\sup_{a\in A}W(a)\,|f(a)-c|\in[0,\infty],
\qquad
\|[f]\|^{\sharp}_{\infty,h;\varphi}:=\inf_{c\in\R}\|f\|_{\infty,h;\varphi;c},
\]
\[
C^{(h)}_\varphi(f;c):=\limsup_{h(a)\to\infty}W(a)\,|f(a)-c|,
\qquad
C^{(h)}_\varphi([f]):=\inf_{c\in\R}C^{(h)}_\varphi(f;c).
\]
For $\varphi=\varphi_p$ we abbreviate $\|\cdot\|_{\infty,h,p;c}$ and
$\|\cdot\|^{\sharp}_{\infty,h,p}$. We define the spaces as the finiteness domains:
\[
C^{c}_{h;\varphi}(A):=\bigl\{f\in C(A):\ \|f\|_{\infty,h;\varphi;c}<\infty\bigr\},
\]
\[
C^{\sharp}_{h;\varphi}(A):=\bigl\{[f]\in C(A)/\R:\ \|[f]\|^{\sharp}_{\infty,h;\varphi}<\infty\bigr\},
\qquad
\Xdec:=\bigl\{[f]\in C^{\sharp}_{h;\varphi}(A):\ C^{(h)}_\varphi([f])=0\bigr\}.
\]
\end{definition}

\begin{remark}[On the definition of the spaces]\label{rem:space-def}
The restriction to the finiteness domain is not cosmetic. On $C(A)/\R$ the sharp
``norm'' takes the value $+\infty$ on most classes and is not a norm; statements such as
completeness or genericity must be made on $C^{\sharp}_{h;\varphi}(A)$.
\end{remark}

\begin{proposition}[Finiteness forces convergence]\label{prop:finite-implies-conv}
Let $\varphi$ be unbounded and $\|f\|_{\infty,h;\varphi;L}=M<\infty$. Then
$|f(a)-L|\le M/W(a)$ and hence $\lim_{a\to\wA}f(a)=L$.
\end{proposition}

\begin{proof}
Immediate from $W(a)=\varphi(h(a))\to\infty$ as $a\to\wA$ (\cref{cor:conv-equiv}).
\end{proof}

\subsection{The reduction theorem}

The following elementary observation determines the whole isotropic theory.

\begin{theorem}[Reduction]\label{thm:reduction}
Let $\varphi\in\Phi_{\mathrm{adm}}$ be unbounded and $f\in C(A)$. The following are
equivalent.
\begin{enumerate}
\item[(i)] $\|[f]\|^{\sharp}_{\infty,h;\varphi}<\infty$;
\item[(ii)] there is $L\in\R$ with $\sup_A W|f-L|<\infty$;
\item[(iii)] $L:=\lim_{a\to\wA}f(a)$ exists in $\R$ and $W\,(f-L)$ is bounded on $A$.
\end{enumerate}
In that case $L$ is unique, and
\[
\|f\|_{\infty,h;\varphi;c}=+\infty\quad\text{for every }c\ne L,
\qquad
\|[f]\|^{\sharp}_{\infty,h;\varphi}=\|f\|_{\infty,h;\varphi;L}.
\]
In particular the infimum over $c$ is attained, at the single point $c=L$.
\end{theorem}

\begin{proof}
(ii)$\Rightarrow$(i) is trivial. (i)$\Rightarrow$(ii): if the infimum is finite there is
$c$ with $\sup_AW|f-c|<\infty$. (ii)$\Leftrightarrow$(iii) is
\cref{prop:finite-implies-conv} together with the definition.

Uniqueness and the last two displays: let $L$ be as in (iii) and let $c\ne L$. Choose
$a_n\to\wA$. Then
\[
W(a_n)|f(a_n)-c|\ \ge\ W(a_n)\bigl(|L-c|-|f(a_n)-L|\bigr)\ \ge\ W(a_n)|L-c|-\sup_AW|f-L|,
\]
and $W(a_n)\to\infty$, so $\|f\|_{\infty,h;\varphi;c}=\infty$. Hence the infimum defining
the sharp norm is over the single admissible value $c=L$.
\end{proof}

\begin{corollary}[Isometry with $\Cb$ and $\Cz$]\label{cor:isometry}
Let $\varphi$ be unbounded. Write $L(f)$ for the limit provided by
\cref{thm:reduction} and define
\[
T:\ C^{\sharp}_{h;\varphi}(A)\longrightarrow \Cb(A),\qquad
T([f]):=W\cdot\bigl(f-L(f)\bigr).
\]
Then $T$ is a well-defined, linear, surjective isometry, and
$T(\Xdec)=\Cz(A)$ isometrically. In particular $C^{\sharp}_{h;\varphi}(A)$ and $\Xdec$
are Banach spaces, and $\Xdec$ is separable whenever $\Cz(A)$ is, e.g.\ when $A$ is
second countable.
\end{corollary}

\begin{proof}
\emph{Well defined and linear.} $L$ is linear on the finiteness domain and satisfies
$L(f+c)=L(f)+c$, so $T([f+c])=W(f+c-L(f)-c)=T([f])$. Linearity of $T$ follows from
linearity of $L$.

\emph{Isometry.} $\|T([f])\|_\infty=\sup_AW|f-L(f)|=\|[f]\|^{\sharp}_{\infty,h;\varphi}$
by \cref{thm:reduction}.

\emph{Surjectivity.} Given $g\in \Cb(A)$ put $f:=g/W\in C(A)$ (recall $W\ge1$ is
continuous). Since $W\to\infty$ at $\wA$ and $g$ is bounded, $f\to0$, so $L(f)=0$ and
$T([f])=g$.

\emph{Decay class.} $C^{(h)}_\varphi([f])=\limsup_{h\to\infty}W|f-L(f)|
=\limsup_{h\to\infty}|T([f])|$, which vanishes exactly when $T([f])\in\Cz(A)$. Conversely
for $g\in\Cz(A)$ the preimage $[g/W]$ lies in $\Xdec$.

Completeness and separability are transported from $\Cb(A)$ and $\Cz(A)$.
\end{proof}

\begin{remark}[Scope of the reduction]\label{rem:scope}
\Cref{thm:reduction} uses only $W\to\infty$ at $\wA$. Three consequences deserve to be
stated plainly.
\begin{enumerate}
\item The theory ``modulo constants'' and the theory ``at fixed $L$'' coincide on the
finiteness domain. There is no genuine quotient phenomenon.
\item Questions of existence, uniqueness, stability and genericity of the minimising
constant are trivial: the minimiser exists, is unique, and equals the limit. The
corresponding statements of the posted version (two-sided contact, generic uniqueness,
upper semicontinuity of the argmin) are either vacuous or false in this setting; see
\cref{app:corrections}.
\item Everything that is a statement about the Banach space
$C^{\sharp}_{h;\varphi}(A)$ or $\Xdec$ alone is a statement about $\Cb(A)$ or $\Cz(A)$.
The weight enters only through the identification $T$, i.e.\ through which \emph{functions}
sit in the space, never through the linear-metric structure.
\end{enumerate}
\end{remark}

\begin{remark}[The bounded case]\label{rem:bounded-phi}
If $\varphi$ is bounded, say $\varphi\equiv1$, then
$\|[f]\|^{\sharp}=\inf_c\sup_A|f-c|=\osc_A(f)$ and one is doing genuine Chebyshev
approximation by constants on a non-compact space. The minimiser may then fail to be
unique, the two-sided contact condition becomes meaningful, and the reduction theorem
fails. This case is genuinely different and is not treated here; see
\cref{q:bounded-phi}.
\end{remark}

\subsection{Splitting the norm and the finiteness criterion}

We record the (correct) form of the ``patching'' identities. The point to be careful about
is that the local and tail infima over $c$ may be attained at different constants, so one
must fix a single $c$ before splitting.

\begin{lemma}[Splitting at a fixed constant]\label{lem:splitting}
For every $c\in\R$ and $R\ge0$,
\[
\|f\|_{\infty,h;\varphi;c}
=\max\Bigl\{\ \sup_{\{h\le R\}}W|f-c|\ ,\ \sup_{\{h>R\}}W|f-c|\ \Bigr\},
\]
and $R\mapsto \sup_{\{h> R\}}W|f-c|$ is non-increasing with
$\lim_{R\to\infty}\sup_{\{h>R\}}W|f-c|=C^{(h)}_\varphi(f;c)$.
\end{lemma}

\begin{proof}
The first identity is the decomposition of a supremum over $A=\{h\le R\}\sqcup\{h>R\}$.
The second is the definition of $\limsup$ along the tail filter.
\end{proof}

\begin{corollary}[Finiteness criterion]\label{cor:finiteness}
Let $\varphi$ be unbounded and $f\in C(A)$. Then
\[
\|[f]\|^{\sharp}_{\infty,h;\varphi}<\infty
\iff
C^{(h)}_\varphi([f])<\infty .
\]
\end{corollary}

\begin{proof}
($\Rightarrow$) is clear. ($\Leftarrow$) If $C^{(h)}_\varphi(f;c)<\infty$ for some $c$
then $W|f-c|$ is bounded on some $\{h>R\}$, hence $|f-c|\le M/W\to0$ there, so $f\to c$.
On the compact set $\{h\le R\}$ the continuous function $W|f-c|$ is bounded. Apply
\cref{lem:splitting}.
\end{proof}

\begin{remark}[On the ``exact patch formula'']
The identity announced in the posted version,
$\|[f]\|^{\sharp}=\inf_{R}\max\{\|[f]\|_{\mathrm{loc},R},T_R([f])\}$ with independent
infima over $c$ in the two terms, is false: taking independent constants can only decrease
the right-hand side, and there is no reason for a single constant to be near-optimal for
both. \Cref{lem:splitting} is the correct substitute, and it is all that is used below.
\end{remark}

\subsection{Comparison of weights}

\begin{proposition}[Global comparison]\label{prop:global-comparison}
Let $\varphi_1\preceq\varphi_2$ in $\Phi_{\mathrm{adm}}$, both unbounded. Then, with $C'$
as in \cref{lem:preceq-global},
\[
\|[f]\|^{\sharp}_{\infty,h;\varphi_1}\ \le\ C'\,\|[f]\|^{\sharp}_{\infty,h;\varphi_2}
\qquad\text{and}\qquad
C^{(h)}_{\varphi_1}(f;L)\ \le\ C\,C^{(h)}_{\varphi_2}(f;L).
\]
In particular $C^{\sharp}_{h;\varphi_2}(A)\hookrightarrow C^{\sharp}_{h;\varphi_1}(A)$
continuously, and $C^{(h)}_{\varphi_2}(f;L)=0$ implies $C^{(h)}_{\varphi_1}(f;L)=0$.
\end{proposition}

\begin{proof}
By \cref{thm:reduction} both sharp norms are computed at the same constant $L=L(f)$;
then $\varphi_1(h)|f-L|\le C'\varphi_2(h)|f-L|$ pointwise. The asymptotic statement uses
only eventual domination.
\end{proof}

\begin{proposition}[Monotonicity in $p$]\label{prop:monotone-p}
For fixed $[f]$, the map $p\mapsto\|[f]\|^{\sharp}_{\infty,h,p}$ is non-decreasing on
$[0,\infty)$.
\end{proposition}

\begin{proof}
$(1+h)^{p}\le(1+h)^{p'}$ pointwise for $p\le p'$ since $1+h\ge1$; apply
\cref{thm:reduction} to evaluate both norms at the same $L$.
\end{proof}

\begin{proposition}[Interpolation for asymptotic constants]\label{prop:interp-asymptotic}
Let $\varphi,\varphi_1,\varphi_2\in\Phi_{\mathrm{adm}}$, $\theta\in[0,1]$, and suppose
$\varphi\le K\varphi_1^{\theta}\varphi_2^{1-\theta}$ on $[0,\infty)$. Assume
$C^{(h)}_{\varphi_1}(f;L)$ and $C^{(h)}_{\varphi_2}(f;L)$ are not respectively $0$ and
$+\infty$ (in either order), so that the right-hand side below is not of the form
$0\cdot\infty$. Then
\[
C^{(h)}_{\varphi}(f;L)\ \le\ K\,\bigl(C^{(h)}_{\varphi_1}(f;L)\bigr)^{\theta}
\bigl(C^{(h)}_{\varphi_2}(f;L)\bigr)^{1-\theta}.
\]
\end{proposition}

\begin{proof}
Pointwise, $\varphi(h)|f-L|\le K(\varphi_1(h)|f-L|)^{\theta}(\varphi_2(h)|f-L|)^{1-\theta}$;
take $\limsup$ and use $\limsup(u_nv_n)\le(\limsup u_n)(\limsup v_n)$, valid for
non-negative sequences except in the indeterminate case $0\cdot\infty$ excluded by
hypothesis (for which $u_n=1/n$, $v_n=n$ is a counterexample).
\end{proof}

\begin{proposition}[Log-convexity in $p$]\label{prop:logconvex}
Let $p_1,p_2>0$, $\theta\in[0,1]$ and $p_\theta=\theta p_1+(1-\theta)p_2$. Then
\[
\|[f]\|^{\sharp}_{\infty,h,p_\theta}\ \le\
\bigl(\|[f]\|^{\sharp}_{\infty,h,p_1}\bigr)^{\theta}
\bigl(\|[f]\|^{\sharp}_{\infty,h,p_2}\bigr)^{1-\theta},
\]
so $p\mapsto\log\|[f]\|^{\sharp}_{\infty,h,p}$ is convex on the interior of its finiteness
interval.
\end{proposition}

\begin{proof}
By \cref{thm:reduction} all three sharp norms are evaluated at the same constant
$L=L(f)$, and pointwise
\[
(1+h)^{p_\theta}|f-L|=\bigl((1+h)^{p_1}|f-L|\bigr)^{\theta}
\bigl((1+h)^{p_2}|f-L|\bigr)^{1-\theta};
\]
take suprema.
\end{proof}

\begin{question}[Log-convexity for bounded weights]\label{q:logconvex}
Does \cref{prop:logconvex} hold for a scale of \emph{bounded} weights, where
\cref{thm:reduction} is unavailable and the three infima over $c$ are attained at
different constants? H\"older gives log-convexity of $c\mapsto\|f\|_{\infty,h;\varphi;c}$
at fixed $c$, but the infimum of a family of log-convex functions need not be log-convex,
and taking $c=c_\ast$ optimal for one exponent gives an inequality in the wrong direction
for the others. The argument of the posted version, which passed from H\"older at fixed $c$
to the infimum over $c$, is invalid here.
\end{question}

\begin{corollary}[Continuity in $p$]\label{cor:cont-p}
For $\varphi_p$, $p>0$, and $[f]$ with $\|[f]\|^{\sharp}_{\infty,h,p_0}<\infty$, the map
$p\mapsto\|[f]\|^{\sharp}_{\infty,h,p}$ is convex in $\log$ hence continuous on the
interior of its finiteness interval, by \cref{prop:logconvex}.
\end{corollary}

\subsection{Elementary operations}

\begin{proposition}[Lipschitz calculus]\label{prop:lipschitz}
Let $\psi:\R\to\R$ be Lipschitz with constant $\Lambda$ and $\varphi$ unbounded. Then
$\|[\psi\circ f]\|^{\sharp}_{\infty,h;\varphi}\le\Lambda\,\|[f]\|^{\sharp}_{\infty,h;\varphi}$,
and $\Xdec$ is stable under $f\mapsto\psi\circ f$.
\end{proposition}

\begin{proof}
$f\to L$ gives $\psi\circ f\to\psi(L)$, so $L(\psi\circ f)=\psi(L(f))$ and
$W|\psi(f)-\psi(L)|\le\Lambda W|f-L|$ pointwise.
\end{proof}

\begin{remark}[No Banach lattice structure]\label{rem:no-lattice}
It was asserted in the posted version that
$\|\,|[f]|\,\|^{\sharp}=\|[f]\|^{\sharp}$ and that $C^{\sharp}_{h;\varphi}(A)$ is a Banach
lattice for the operations induced by $\max$ and $\min$. Only the inequality
$\|\,|[f]|\,\|^{\sharp}\le\|[f]\|^{\sharp}$ holds (a case of \cref{prop:lipschitz}), and
it can be strict: on $A=[0,\infty)$ with $h(x)=x$, $\varphi(t)=1+t$ and
$f(x)=1-2e^{-x}$ one has $L(f)=1$, $\|[f]\|^{\sharp}\ge|f(0)-1|=2$, whereas
$|f|(0)=1$ and $L(|f|)=1$, so the contribution of the point $0$ to
$\|\,|[f]|\,\|^{\sharp}$ vanishes. Moreover the order on $C(A)/\R$ induced by pointwise
order is not well defined (adding constants does not preserve it), so the lattice
statement is not even meaningful as written.
\end{remark}

\begin{theorem}[Weighted Schur test]\label{thm:schur}
Let $(A,\mu)$ be a Radon measure space, $K\ge0$ measurable with
$\int_AK(x,y)\,d\mu(y)=1$ for every $x$ (so that $T_K$ preserves constants), and set
$(T_Kf)(x)=\int_AK(x,y)f(y)\,d\mu(y)$, assumed to define a map $C(A)\to C(A)$ (e.g.\ if
$K(\cdot,y)$ is continuous for each $y$ and a local domination hypothesis allows passing to
the limit under the integral sign). If
\[
C_1:=\sup_{x\in A}\int_A K(x,y)\,\frac{W(x)}{W(y)}\,d\mu(y)<\infty,
\]
then $T_K$ maps $C^{\sharp}_{h;\varphi}(A)$ boundedly into itself with
$\|[T_Kf]\|^{\sharp}\le C_1\|[f]\|^{\sharp}$, and $L(T_Kf)=L(f)$.
\end{theorem}

\begin{proof}
Let $L=L(f)$. Using $\int K(x,\cdot)d\mu=1$,
\[
|T_Kf(x)-L|=\Bigl|\int K(x,y)(f(y)-L)\,d\mu(y)\Bigr|
\le \|[f]\|^{\sharp}\int K(x,y)W(y)^{-1}d\mu(y),
\]
so $W(x)|T_Kf(x)-L|\le C_1\|[f]\|^{\sharp}$. Since $W\to\infty$, this forces
$T_Kf\to L$, i.e.\ $L(T_Kf)=L$; then \cref{thm:reduction} identifies the sharp norm.
\end{proof}

\begin{remark}
The posted version imposed instead $\sup_y\int K(x,y)d\mu(x)\le C_2$, which is never used,
and the normalisation $\int K(x,y)d\mu(x)=1$, which is the wrong variable: without
$\int K(x,\cdot)\,d\mu=1$ the operator does not preserve constants and does not descend
to the quotient.
\end{remark}

\begin{lemma}[Local uniform convergence plus tight tails]\label{lem:dini}
Let $c\in\R$, and let $f_n,f\in C^{c}_{h;\varphi}(A)$ satisfy
\begin{enumerate}
\item[(i)] $\sup_{\{h\le R\}}W|f_n-f|\to0$ for every $R<\infty$;
\item[(ii)] $\lim_{R\to\infty}\ \sup_{n}\ \sup_{\{h>R\}}W|f_n-f|=0$.
\end{enumerate}
Then $\sup_AW|f_n-f|\to0$; in particular $\|[f_n]-[f]\|^{\sharp}\to0$.
\end{lemma}

\begin{proof}
Given $\varepsilon>0$ pick $R$ with $\sup_n\sup_{\{h>R\}}W|f_n-f|<\varepsilon$ by (ii),
then $N$ with $\sup_{\{h\le R\}}W|f_n-f|<\varepsilon$ for $n\ge N$ by (i). Conclude with
\cref{lem:splitting}.
\end{proof}

\begin{proposition}[Density of truncations]\label{prop:truncation}
Let $\varphi$ be unbounded. For every $[f]\in\Xdec$ and $\varepsilon>0$ there exist a
compact $K\subset A$, $c\in\R$ and $g\in C(A)$ with $g=f$ on $K$, $g\equiv c$ off a compact
set, and $\|[f]-[g]\|^{\sharp}_{\infty,h;\varphi}<\varepsilon$.
\end{proposition}

\begin{proof}
Via \cref{cor:isometry} this is the density of $C_c(A)$ in $\Cz(A)$: choose
$u\in C_c(A)$ with $\|T([f])-u\|_\infty<\varepsilon$ and $u=T([f])$ on a prescribed
compact set (multiply $T([f])$ by a Urysohn cut-off equal to $1$ on $\{h\le R_0\}$ and
supported in $\{h\le R_1\}$; the error is bounded by $\sup_{\{h>R_0\}}|T([f])|<\varepsilon$
for $R_0$ large, since $T([f])\in\Cz(A)$). Then $g:=L(f)+u/W$ works, and
$\|[f]-[g]\|^{\sharp}=\|T([f])-u\|_\infty$.
\end{proof}

\subsection{Coarse robustness: dilation along the image of $h$}

Two invariance properties may be asked of the asymptotic constants. Fix
$\varphi\in\Phi_{\mathrm{adm}}$ and write, for a continuous proper $h'$ with $h'\sim h$,
\begin{itemize}
\item[(INV)] for every such $h'$ there is $C\ge1$ with
$C^{-1}C^{(h)}_\varphi(f;L)\le C^{(h')}_\varphi(f;L)\le C\,C^{(h)}_\varphi(f;L)$ for all
$f\in C(A)$, $L\in\R$;
\item[(VAN)] for every such $h'$ and all $f,L$:
$C^{(h)}_\varphi(f;L)=0\Rightarrow C^{(h')}_\varphi(f;L)=0$.
\end{itemize}
Clearly (INV) implies (VAN). It turns out that both are governed by how much $\varphi$ can
be dilated \emph{at the scales that $h$ actually visits}, and that they are not equivalent.

\begin{definition}[Dilation functionals]\label{def:bounded-dilation}
For $\lambda\ge1$ set
\[
D_h(\lambda):=\sup_{t\in h(A)}\frac{\varphi(\lambda t)}{\varphi(t)},
\qquad
E_h(\lambda):=\sup_{t\in h(A)}\frac{\varphi(t)}{\varphi(t/\lambda)},
\]
\[
D(\lambda):=\sup_{t\ge0}\frac{\varphi(\lambda t)}{\varphi(t)} .
\]
All three lie in $[1,\infty]$ and are non-decreasing in $\lambda$. We say $\varphi$ has
\emph{globally bounded dilation}, or is \emph{$O$-regularly varying}
\cite{BinghamGoldieTeugels1987}, if $D(\lambda)<\infty$ for every $\lambda\ge1$; note
$D_h\le D$ and $E_h\le D$, the latter by the substitution $s=t/\lambda$. The polynomial
weights $\varphi_p$ ($D(\lambda)=\lambda^{p}$), the log-polynomial weights
$\varphi_{p,q}$ and every bounded $\varphi$ have globally bounded dilation; the
exponential weights $\varphi_\alpha(t)=e^{\alpha t}$, $\alpha>0$, do not, since
$\varphi_\alpha(\lambda t)/\varphi_\alpha(t)=e^{\alpha(\lambda-1)t}\to\infty$.
\end{definition}

\begin{theorem}[Characterisation of coarse-affine invariance]\label{thm:coarse-stability}
Let $h$ be a continuous proper exhaustion of $A$ and $\varphi\in\Phi_{\mathrm{adm}}$.
\begin{enumerate}
\item[(a)] \textup{(VAN)} holds if and only if $D_h(\lambda)<\infty$ for every $\lambda\ge1$.
\item[(b)] \textup{(INV)} holds if and only if $D_h(\lambda)<\infty$ and
$E_h(\lambda)<\infty$ for every $\lambda\ge1$.
\item[(c)] If $h(A)\supseteq[t_0,\infty)$ for some $t_0<\infty$ --- in particular if $A$ is
connected, since then $h(A)$ is an unbounded interval --- then \textup{(INV)},
\textup{(VAN)} and globally bounded dilation of $\varphi$ are all equivalent.
\end{enumerate}
\end{theorem}

\begin{proof}
Throughout, $h\sim h'$ implies that $h\to\infty$ exactly when $h'\to\infty$, so the two
$\limsup$'s are taken along the same filter. Write $W=\varphi\circ h$, $W'=\varphi\circ h'$.

\emph{(a), sufficiency.} Let $a_1h-b_1\le h'\le a_2h+b_2$. By (A2) and monotonicity,
\[
W'=\varphi(h')\le\varphi(a_2h+b_2)\le K\varphi(b_2)\,\varphi(a_2h)
\le K\varphi(b_2)\,D_h(\max(1,a_2))\,W ,
\]
using $D_h$ at the points $h(a)\in h(A)$ (for $a_2<1$ monotonicity suffices). Multiplying
by $|f-L|$ and taking $\limsup$ gives $C^{(h')}_\varphi\le C\,C^{(h)}_\varphi$, hence (VAN).

\emph{(a), necessity.} Suppose $D_h(\lambda_0)=\infty$ for some $\lambda_0\ge1$. Choose
$t_n\in h(A)$, which after passing to a subsequence we may take strictly increasing with
$t_n\to\infty$, such that $r(t_n):=\varphi(\lambda_0t_n)/\varphi(t_n)\ge n^{2}$; this is
possible because $\varphi(\lambda_0\,\cdot)/\varphi$ is bounded on compact intervals
($\varphi$ is continuous and $\ge1$), so the supremum can only be infinite along a sequence
tending to infinity. Put $h':=\lambda_0h$, a continuous proper exhaustion with $h'\sim h$.
Pick a continuous non-decreasing $v:[0,\infty)\to[1,\infty)$ with $v(t_n)=n$ and
$v(t)\to\infty$ (linear interpolation between the $t_n$, which is legitimate because they
are strictly increasing), and set
\[
f:=L+\frac{1}{\varphi(h)\,v(h)}\ \in C(A).
\]
Then $\varphi(h)|f-L|=1/v(h)\to0$, so $C^{(h)}_\varphi(f;L)=0$. Since $t_n\in h(A)$ there
is $a_n\in A$ with $h(a_n)=t_n$, and
\[
\varphi(h'(a_n))\,|f(a_n)-L|=\frac{\varphi(\lambda_0t_n)}{\varphi(t_n)\,v(t_n)}
=\frac{r(t_n)}{n}\ \ge\ n\longrightarrow\infty,
\]
so $C^{(h')}_\varphi(f;L)=\infty$ and (VAN) fails. Note that this argument uses
$t_n\in h(A)$ in an essential way: it is exactly the step that requires the bad scales to
be visible to $h$.

\emph{(b), sufficiency.} The upper bound $W'\le CW$ is as in (a). For the lower bound, let
$\lambda:=\max(1,a_1^{-1})$. On $\{h\ge b_1/a_1\}$ we have $a_1h-b_1\ge0$ and, by (A2),
$\varphi(a_1h)=\varphi\bigl((a_1h-b_1)+b_1\bigr)\le K\varphi(b_1)\varphi(a_1h-b_1)$, while
$W=\varphi(h)\le E_h(\lambda)\,\varphi(h/\lambda)\le E_h(\lambda)\varphi(a_1h)$ if
$a_1\le1$ (and $W\le\varphi(a_1h)$ if $a_1\ge1$, by monotonicity). Combining,
$W\le E_h(\lambda)K\varphi(b_1)\varphi(a_1h-b_1)\le E_h(\lambda)K\varphi(b_1)W'$. On the
complementary set $\{h<b_1/a_1\}$, which is compact, both $W$ and $W'$ are bounded above
and below by positive constants. Hence $W\asymp W'$ on $A$, and (INV) follows.

\emph{(b), necessity.} For $\mu>0$ the function $h':=\mu h$ is a continuous proper
exhaustion with $h'\sim h$. Applying (INV) to $f:=L+1/\varphi(h)$ gives
$C^{(h)}_\varphi(f;L)=1$ and
$C^{(h')}_\varphi(f;L)=\limsup\varphi(\mu h)/\varphi(h)$, so this $\limsup$ is finite and
non-zero; running the same argument with $f:=L+1/\varphi(\mu h)$ gives the reciprocal
bound. Since $\varphi(\mu\,\cdot)/\varphi$ is bounded above and below by positive constants
on compacta, we conclude $\sup_{t\in h(A)}\varphi(\mu t)/\varphi(t)<\infty$ and
$\sup_{t\in h(A)}\varphi(t)/\varphi(\mu t)<\infty$ for every $\mu>0$, i.e.\ $D_h$ and
$E_h$ are finite.

\emph{(c).} If $h(A)\supseteq[t_0,\infty)$ then $D_h(\lambda)<\infty$ forces
$D(\lambda)<\infty$, since $\sup_{t\le t_0}\varphi(\lambda t)/\varphi(t)\le\varphi(\lambda t_0)<\infty$.
Conversely $D(\lambda)<\infty$ gives both $D_h(\lambda)\le D(\lambda)$ and
$E_h(\lambda)\le D(\lambda)$. Now apply (a) and (b).
\end{proof}

\begin{corollary}[Polynomial and log-polynomial weights]\label{cor:poly-invariance}
For $\varphi=\varphi_p$ or $\varphi=\varphi_{p,q}$, and any continuous proper $h'\sim h$,
the constants $C^{(h)}_\varphi(f;L)$ and $C^{(h')}_\varphi(f;L)$ are comparable, with
constants depending only on $p,q$ and the coarse-affine data. In particular the decay class
$O(h^{-p})$ depends only on the coarse-affine class of $h$, for every $A$.
\end{corollary}

\begin{proof}
These weights have globally bounded dilation, so \cref{def:bounded-dilation} gives
$D_h,E_h<\infty$ and \cref{thm:coarse-stability}(b) applies.
\end{proof}

\begin{example}[Exponential weights, connected case]\label{ex:exp-not-invariant}
Take $A=[0,\infty)$ (connected), $h(x)=x$, $h'(x)=2x$, $\varphi(t)=e^{t}$ and
$f-L=e^{-x}$. Then
\[
C^{(h)}_\varphi(f;L)=\limsup_{x\to\infty}e^{x}e^{-x}=1,
\qquad
C^{(h')}_\varphi(f;L)=\limsup_{x\to\infty}e^{2x}e^{-x}=+\infty .
\]
With $f-L=e^{-x}/(1+x)$ one gets $C^{(h)}_\varphi=0$ and $C^{(h')}_\varphi=\infty$, so
(VAN) fails as well. Here $h(A)=[0,\infty)$, so this is forced by
\cref{thm:coarse-stability}(c).
\end{example}

\subsubsection{Connectedness is not cosmetic: two examples on $\N$}

The equivalences of \cref{thm:coarse-stability}(c) genuinely require $h$ to see all large
scales. We record two counterexamples on the discrete space $A=\N$ of
\cref{ex:instances}, built from the same construction.

\begin{lemma}[Staircase weights]\label{lem:staircase}
Let $0<\alpha_1<\beta_1<\alpha_2<\beta_2<\cdots$ with $\alpha_{n+1}/\beta_n\to\infty$, and
let $a_n>0$ satisfy, for all large $n$,
\[
\text{(C1)}\quad a_n\le s_{n-1}:=\sum_{m<n}a_m,
\qquad
\text{(C2)}\quad \frac{a_n\,\beta_{n-1}}{\beta_n-\alpha_n}\le c_0<\infty .
\]
Let $\psi:[0,\infty)\to[0,\infty)$ be continuous, non-decreasing, equal to $0$ on
$[0,\alpha_1]$, affine with slope $\sigma_n:=a_n/(\beta_n-\alpha_n)$ on $[\alpha_n,\beta_n]$
and constant on $[\beta_n,\alpha_{n+1}]$. Then $\varphi:=e^{\psi}\in\Phi_{\mathrm{adm}}$,
and $\varphi$ is unbounded.
\end{lemma}

\begin{proof}
(A1) is clear, with $\varphi(0)=1$, and $\psi(\beta_n)=s_n\to\infty$. For (A2) we must find
$K$ with $\psi(r+s)\le\log K+\psi(r)+\psi(s)$; by symmetry assume $r\le s$, so it suffices
to bound the increment $\Delta:=\psi(s+r)-\psi(s)$ by $\log K+\psi(r)$. Since
$s+r\le2s$ and $\alpha_{n+1}/\beta_n\to\infty$, for $s$ large the interval $[s,s+r]$ meets
at most one rise interval $[\alpha_n,\beta_n]$, and then $\Delta\le\sigma_n r$. If
$r\ge\beta_{n-1}$ then $\psi(r)\ge s_{n-1}\ge a_n\ge\Delta$ by (C1). If $r<\beta_{n-1}$
then $\Delta\le\sigma_n\beta_{n-1}=a_n\beta_{n-1}/(\beta_n-\alpha_n)\le c_0$ by (C2). The
finitely many remaining values of $n$, and the bounded range of $s$, contribute a further
constant.
\end{proof}

\begin{example}[(INV) holds although $\varphi$ has unbounded global dilation]
\label{ex:N-inv-strict}
Let $A=\N$ with $h(k)=k!$, a continuous (all functions on a discrete space are) proper
exhaustion, and $h(A)=\{k!:k\in\N\}$. Apply \cref{lem:staircase} with
\[
\alpha_n=n!\,n^{1/4},\qquad \beta_n=n!\,\sqrt n,\qquad a_n=n .
\]
Then $\alpha_{n+1}/\beta_n=(n+1)^{5/4}/\sqrt n\to\infty$, (C1) holds since
$s_{n-1}\sim n^{2}/2\ge n$, and (C2) holds since
$a_n\beta_{n-1}/(\beta_n-\alpha_n)=n(n-1)!\sqrt{n-1}/\bigl(n!(\sqrt n-n^{1/4})\bigr)\to1$.
Then:
\begin{itemize}
\item \emph{Global dilation is unbounded:} at $t=\alpha_n$ one has $[t,2t]\subset
[\alpha_n,\beta_n]$ for $n\ge16$, so
$\log\bigl(\varphi(2t)/\varphi(t)\bigr)=\sigma_n\alpha_n=n^{5/4}/(\sqrt n-n^{1/4})\to\infty$,
i.e.\ $D(2)=\infty$.
\item \emph{Yet \textup{(INV)} holds:} if $h'\sim h$ then $h'(k)\in[a_1k!-b_1,a_2k!+b_2]$,
and for $k$ large this window lies strictly between $\beta_{k-1}=(k-1)!\sqrt{k-1}$
(which is $o(k!)$) and $\alpha_k=k!\,k^{1/4}$ (which is $\gg a_2k!$), i.e.\ inside a
plateau of $\psi$ that also contains $h(k)=k!$. Hence $\varphi(h'(k))=\varphi(h(k))$ for
all large $k$, so $C^{(h')}_\varphi(f;L)=C^{(h)}_\varphi(f;L)$ for every $f$ and $L$.
\end{itemize}
Consistently with \cref{thm:coarse-stability}(b), $D_h$ and $E_h$ are finite here: for
$\lambda\ge1$ neither $[k!,\lambda k!]$ nor $[k!/\lambda,k!]$ meets a rise interval, once
$k$ is large.
\end{example}

\begin{example}[(VAN) holds but (INV) fails]\label{ex:N-van-strict}
Again $A=\N$, $h(k)=k!$, and \cref{lem:staircase} with
\[
\alpha_n=n!/2,\qquad \beta_n=n!,\qquad a_n=n,
\]
for which $\alpha_{n+1}/\beta_n=(n+1)/2\to\infty$, (C1) holds, and
$a_n\beta_{n-1}/(\beta_n-\alpha_n)=n(n-1)!/(n!/2)=2$, so (C2) holds. Here
$\psi(k!)=s_k$ and $\psi(k!/2)=s_{k-1}$.
\begin{itemize}
\item \emph{$D_h$ is finite:} for $\lambda\ge1$, the interval $[k!,\lambda k!]$ contains no
rise interval once $k+1>2\lambda$, so $D_h(\lambda)<\infty$; by
\cref{thm:coarse-stability}(a), \textup{(VAN)} holds.
\item \emph{$E_h$ is infinite and \textup{(INV)} fails:}
$E_h(2)\ge\varphi(k!)/\varphi(k!/2)=e^{a_k}=e^{k}\to\infty$. Explicitly, take
$h':=h/2\sim h$ and $f(k):=e^{-s_k}$, $L=0$. Then
\[
C^{(h)}_\varphi(f;0)=\limsup_k e^{s_k}e^{-s_k}=1,
\qquad
C^{(h')}_\varphi(f;0)=\limsup_k e^{s_{k-1}}e^{-s_k}=\limsup_k e^{-k}=0 .
\]
\end{itemize}
So the two invariance properties are genuinely distinct, and neither is implied by
finiteness of $D_h$ alone together with the other.
\end{example}

\begin{remark}[Two false statements corrected during revision]\label{rem:invariance-history}
An intermediate draft asserted, for every weight with ``power-dilation control''
($\varphi(\lambda t)\le C(\lambda)\varphi(t)^{\kappa(\lambda)}$, $\kappa\ge1$ allowed, a
class it claimed included the exponential weights), the estimate
$C^{(h')}_\varphi\le C_1(C^{(h)}_\varphi)^{\kappa_1}$. This is false by
\cref{ex:exp-not-invariant}. The invalid step was
``$W|f-L|\le M+\varepsilon$, whence $W^{\kappa}|f-L|\le(M+\varepsilon)^{\kappa}$ using
$W\ge1$'': in fact
$W^{\kappa}|f-L|=(W|f-L|)^{\kappa}\,|f-L|^{1-\kappa}$, and for $\kappa>1$ the second factor
diverges precisely because $f\to L$; the claimed inequality is equivalent to $|f-L|\ge1$.

A later draft replaced this by an equivalence between globally bounded dilation and
(INV)/(VAN), but its proof of necessity chose points $a_n$ with $h(a_n)=t_n$, which
presupposes that $h$ attains all large values. That holds when $A$ is connected and fails
in general; \cref{ex:N-inv-strict} shows the equivalence itself is false on $A=\N$. It is
also not enough to replace the global dilation by the dilation along $h(A)$:
\cref{ex:N-van-strict} has $D_h$ finite while (INV) fails. \Cref{thm:coarse-stability}
separates the two properties and characterises each; \cref{cor:poly-invariance}, which is
what the rest of the paper uses, is unaffected by either error.
\end{remark}

\begin{corollary}\label{cor:coarse-iso}
If $h\sim h'$ and $\varphi=\varphi_p$, the identity map induces a linear isomorphism
$C^{\sharp}_{h;\varphi_p}(A)\to C^{\sharp}_{h';\varphi_p}(A)$ with equivalent norms
(\cref{lem:equiv-affine}).
\end{corollary}

\begin{remark}
By \cref{ex:no-coarse-affine}, \cref{cor:poly-invariance} is a statement about a
\emph{fixed} coarse-affine class, not an absolute invariance: distinct proper exhaustions
of the same space may lie in distinct classes and then define distinct decay classes. The
corresponding claim of the posted version, that $O(h^{-p})$ is independent of any reasonable
choice of exhaustion, is false. Thus exhaustions are classified up to coarse-affine
equivalence and weights up to dilation ---  but the relevant dilation is the one visible
along $h(A)$, and even that splits into two distinct conditions
(\cref{thm:coarse-stability}).
\end{remark}

\section{Duality and compactness}
\label{sec:duality}

Throughout this section $\varphi\in\Phi_{\mathrm{adm}}$ is unbounded, $W=\varphi\circ h$,
and $T$ is the isometry of \cref{cor:isometry}. We assume in addition that $A$ is second
countable, so that $\Cz(A)$ is separable and Riesz representation applies in the usual
form.

\subsection{The dual of the decay class}

\begin{definition}[Weighted zero-mass measures]\label{def:M0}
Let $\Msp_0^{W^{-1}}(A)$ be the space of finite signed regular Borel measures $\nu$ on $A$
with
\[
\nu(A)=0,\qquad \|\nu\|_{W^{-1}}:=\int_AW^{-1}\,d|\nu|<\infty .
\]
(The integrability is automatic since $W\ge1$ and $|\nu|(A)<\infty$; the condition
$\nu(A)=0$ is what makes $\int f\,d\nu$ depend only on $[f]$.)
\end{definition}

\begin{proposition}[Canonical embedding]\label{prop:embedding}
For $\nu\in\Msp_0^{W^{-1}}(A)$ the formula $\Lambda_\nu([f]):=\int_Af\,d\nu$ defines a
bounded linear functional on $\Xdec$, and under the Riesz identification
$\Xdec^{\ast}\cong\Msp(A)$ induced by $T$ it corresponds to the measure
$\sigma_\nu:=W^{-1}\nu$. Moreover $\|\Lambda_\nu\|_{\Xdec^{\ast}}=\|\nu\|_{W^{-1}}$.
\end{proposition}

\begin{proof}
Since $\nu(A)=0$, $\int f\,d\nu=\int(f-L(f))\,d\nu=\int W(f-L(f))\,W^{-1}d\nu
=\int T([f])\,d\sigma_\nu$. As $T:\Xdec\to\Cz(A)$ is a surjective isometry and
$\Cz(A)^{\ast}=\Msp(A)$ isometrically (Riesz), the functional $\Lambda_\nu$ has norm
$|\sigma_\nu|(A)=\int W^{-1}d|\nu|=\|\nu\|_{W^{-1}}$.
\end{proof}

\begin{theorem}[Duality: isometric, dense, not surjective]\label{thm:duality}
The map $\nu\mapsto\Lambda_\nu$ is an isometric linear embedding of
$\bigl(\Msp_0^{W^{-1}}(A),\|\cdot\|_{W^{-1}}\bigr)$ into $\Xdec^{\ast}$. Its range is,
under the identification of \cref{prop:embedding},
\[
\mathcal{R}=\Bigl\{\sigma\in\Msp(A):\ \int_AW\,d|\sigma|<\infty\ \text{ and }\ \int_AW\,d\sigma=0\Bigr\},
\]
which is norm-dense in $\Msp(A)$ but proper. In particular $\Xdec^{\ast}$ is
\emph{not} isometrically isomorphic to $\Msp_0^{W^{-1}}(A)$; it is isometrically
isomorphic to $\Msp(A)$.
\end{theorem}

\begin{proof}
Isometry is \cref{prop:embedding}. The correspondence $\nu\leftrightarrow\sigma=W^{-1}\nu$
is a bijection between $\{\nu\ \text{finite signed},\ \nu(A)=0\}$ and
$\{\sigma:\int W\,d|\sigma|<\infty,\ \int W\,d\sigma=0\}$: indeed $\nu=W\sigma$ is finite
iff $\int W\,d|\sigma|<\infty$, and $\nu(A)=\int W\,d\sigma$. This identifies
$\mathcal{R}$.

\emph{Properness.} Fix $a_0\in A$ and take $\sigma=\delta_{a_0}$. Then
$\int W\,d\sigma=W(a_0)\ge1\ne0$, so $\delta_{a_0}\notin\mathcal{R}$.

\emph{Density.} Let $\sigma\in\Msp(A)$ and $\varepsilon>0$. Compactly supported measures
are norm-dense in $\Msp(A)$ and satisfy $\int W\,d|\sigma|<\infty$ ($W$ is bounded on
compacta), so we may assume $\supp\sigma$ compact. Set $m:=\int W\,d\sigma$. Choose
$b\in A$ with $W(b)>|m|/\varepsilon$ (possible since $W$ is unbounded) and put
$\sigma':=\sigma-\dfrac{m}{W(b)}\delta_b$. Then $\int W\,d\sigma'=m-m=0$,
$\int W\,d|\sigma'|<\infty$, so $\sigma'\in\mathcal{R}$, and
$\|\sigma-\sigma'\|_{\Msp}=|m|/W(b)<\varepsilon$.
\end{proof}

\begin{corollary}[Norming by measures, and the failure of extremals]
\label{cor:norming}
For every $[f]\in\Xdec$,
\[
\|[f]\|^{\sharp}_{\infty,h;\varphi}
=\sup\Bigl\{\Bigl|\int_Af\,d\nu\Bigr|:\ \nu\in\Msp_0^{W^{-1}}(A),\ \|\nu\|_{W^{-1}}\le1\Bigr\},
\]
but the supremum need not be attained.
\end{corollary}

\begin{proof}
By \cref{thm:duality} the unit ball of $\Msp_0^{W^{-1}}(A)$ is norm-dense in the unit ball
of $\Xdec^{\ast}$, and $\|[f]\|^{\sharp}=\sup_{\|\Lambda\|\le1}|\Lambda([f])|$ by
Hahn--Banach. Non-attainment is \cref{ex:no-extremal}.
\end{proof}

\begin{example}[No extremal measure]\label{ex:no-extremal}
Take $A=[0,\infty)$, $h(x)=x$, $\varphi(t)=1+t$, so $W(x)=1+x$, and
\[
f(x):=(1+x)^{-2},\qquad\text{so}\qquad g:=T([f])=(1+x)^{-1}\in\Cz(A),
\]
whence $L(f)=0$, $[f]\in\Xdec$ and $\|[f]\|^{\sharp}=\|g\|_\infty=1$, attained only at
$x=0$. Let $\nu\in\Msp_0^{W^{-1}}(A)$ with $\|\nu\|_{W^{-1}}\le1$ and put
$\sigma:=W^{-1}\nu$, so $|\sigma|(A)\le1$ and $\int W\,d\sigma=\nu(A)=0$. Then
$\int f\,d\nu=\int g\,d\sigma$, and
\[
\Bigl|\int g\,d\sigma\Bigr|\le\int|g|\,d|\sigma|\le|\sigma|(A)\le1,
\]
with equality throughout only if $|\sigma|$ is carried by $\{|g|=1\}=\{0\}$ and
$|\sigma|(A)=1$, i.e.\ $\sigma=\pm\delta_0$. But then $\int W\,d\sigma=\pm W(0)=\pm1\ne0$,
contradicting $\nu(A)=0$. Hence the supremum $1$ is approached but never attained.
\end{example}

\begin{proposition}[Exact criterion for the existence of a norming measure]
\label{prop:extremal-iff}
Let $[f]\in\Xdec$, $[f]\ne0$, and put $g:=T([f])\in\Cz(A)$, $M:=\|g\|_\infty>0$. The
supremum in \cref{cor:norming} is attained if and only if
\[
\max_A g=M\quad\text{and}\quad\min_A g=-M,
\]
i.e.\ if and only if $g$ attains both $+M$ and $-M$.
\end{proposition}

\begin{proof}
Recall from \cref{prop:embedding} that admissible competitors correspond to
$\sigma\in\Msp(A)$ with $|\sigma|(A)\le1$ and $\int W\,d\sigma=0$, the value being
$\int g\,d\sigma$.

\emph{Sufficiency.} Suppose $g(a)=M$ and $g(b)=-M$. Put
\[
\sigma:=t\,\delta_a-(1-t)\,\delta_b,\qquad t:=\frac{W(b)}{W(a)+W(b)}\in(0,1).
\]
Then $|\sigma|(A)=1$ and $\int W\,d\sigma=tW(a)-(1-t)W(b)=0$, while
$\int g\,d\sigma=tM+(1-t)M=M$.

\emph{Necessity.} Suppose $\sigma$ is admissible with $\bigl|\int g\,d\sigma\bigr|=M$;
replacing $\sigma$ by $-\sigma$ (which preserves both constraints) we may assume
$\int g\,d\sigma=M$. Writing $\sigma=\sigma^{+}-\sigma^{-}$ for the Jordan decomposition,
\[
M=\int g\,d\sigma^{+}-\int g\,d\sigma^{-}\ \le\ M\,\sigma^{+}(A)+M\,\sigma^{-}(A)
=M\,|\sigma|(A)\ \le\ M .
\]
Equality throughout forces $|\sigma|(A)=1$, $\sigma^{+}$ carried by $\{g=M\}$ and
$\sigma^{-}$ carried by $\{g=-M\}$. If $\sigma^{-}=0$ then $\sigma=\sigma^{+}\ge0$ with
$\sigma(A)=1$, so $\int W\,d\sigma\ge\sigma(A)=1>0$ since $W\ge1$, contradicting
$\int W\,d\sigma=0$; symmetrically $\sigma^{+}\ne0$. Hence both $\{g=M\}$ and $\{g=-M\}$
are non-empty.
\end{proof}

\begin{remark}
Since $g\in\Cz(A)$, it attains its maximum whenever that maximum is positive and its
minimum whenever that minimum is negative. \Cref{prop:extremal-iff} therefore says that a
norming measure exists exactly when $g$ is ``balanced'', $\max g=-\min g$. In
\cref{ex:no-extremal}, $g=(1+x)^{-1}>0$ has $\min g=0\ne-M$, which is precisely why the
supremum is not attained. Non-attainment is thus the generic situation, and
\cref{prop:extremal-iff} delimits the exception exactly.
\end{remark}

\begin{remark}[Retracted statements]\label{rem:retracted-duality}
The posted version asserted (a) a full duality $\Msp_0^{W^{-1}}(A)\cong\Xdec^{\ast}$,
(b) sequential weak-$*$ compactness of the unit ball of $\Msp_0^{W^{-1}}(A)$, and (c)
existence of dual extremals and a description of the subdifferential of the sharp norm
supported on a two-sided contact set with masses $\pm\tfrac12$. All three are false:
(a) by \cref{thm:duality}, (b) because a norm-dense proper subspace has a unit ball which
is weak-$*$ dense but not weak-$*$ closed in the ball of the dual, and (c) by
\cref{ex:no-extremal} together with the observation that the contact set need not have
points of both signs (take $f\ge L$). The arguments given there were also mutually
circular: (a) invoked (b) and (b) invoked (a).
\end{remark}

\begin{remark}[What is true about weak-$*$ compactness]
Since $\Xdec\cong\Cz(A)$ is separable, the closed unit ball of $\Xdec^{\ast}\cong\Msp(A)$
is weak-$*$ sequentially compact by Banach--Alaoglu and metrisability. This is the
statement to use in practice; it is about $\Msp(A)$, not about $\Msp_0^{W^{-1}}(A)$.
\end{remark}

\subsection{Compactness}

\begin{theorem}[Weighted Arzel\`a--Ascoli]\label{thm:AA}
Let $\mathcal{F}\subset\Xdec$ and put $g_{[f]}:=T([f])=W(f-L(f))$. Suppose
\begin{enumerate}
\item[(i)] $\sup_{[f]\in\mathcal{F}}\|[f]\|^{\sharp}_{\infty,h;\varphi}<\infty$;
\item[(ii)] $\{g_{[f]}:[f]\in\mathcal{F}\}$ is equicontinuous on every compact subset of $A$;
\item[(iii)] tight tails: $\displaystyle\lim_{R\to\infty}\ \sup_{[f]\in\mathcal{F}}\ \sup_{\{h>R\}}W\,|f-L(f)|=0$.
\end{enumerate}
Then $\mathcal{F}$ is relatively compact in $\bigl(\Xdec,\|\cdot\|^{\sharp}\bigr)$.
\end{theorem}

\begin{proof}
By \cref{cor:isometry} it suffices to prove that $G:=T(\mathcal{F})\subset\Cz(A)$ is
relatively compact for the sup-norm. Let $\varepsilon>0$. By (iii) choose $R$ with
$|g|<\varepsilon/3$ on $\{h>R\}$ for all $g\in G$. On the compact set $K_R=\{h\le R\}$ the
family $G|_{K_R}$ is bounded by (i) and equicontinuous by (ii), hence relatively compact in
$C(K_R)$ by the classical Arzel\`a--Ascoli theorem; pick $g_1,\dots,g_N\in G$ whose
restrictions form an $\varepsilon/3$-net for $G|_{K_R}$. For $g\in G$ choose $j$ with
$\sup_{K_R}|g-g_j|<\varepsilon/3$; then
$\sup_A|g-g_j|\le\max\{\varepsilon/3,\ 2\varepsilon/3\}<\varepsilon$. Thus $G$ is totally
bounded in $\Cz(A)$, which is complete.
\end{proof}

\begin{remark}
The hypothesis of equicontinuity must be imposed on $g_{[f]}=W(f-L(f))$, not on $f$: it is
$g_{[f]}$ that lives in the space in which compactness is asserted. If $W$ is locally
Lipschitz and the $L(f)$ range in a bounded set, equicontinuity of $\{f\}$ on compacta
implies (ii).
\end{remark}

\begin{theorem}[Compact inclusion]\label{thm:compact-inclusion}
Let $\varphi_1\preceq\varphi_2$ in $\Phi_{\mathrm{adm}}$ be unbounded with
$\varphi_1(t)/\varphi_2(t)\to0$. Let $\mathcal{F}$ be bounded in
$\|\cdot\|^{\sharp}_{\infty,h;\varphi_2}$ and such that
$\{W_2(f-L(f))\}$ is equicontinuous on compacta, where $W_j=\varphi_j\circ h$. Then
$\mathcal{F}$ is relatively compact in $C^{\sharp}_{h;\varphi_1}(A)$.
\end{theorem}

\begin{proof}
Let $M:=\sup_{\mathcal F}\|[f]\|^{\sharp}_{\infty,h;\varphi_2}$. By \cref{thm:reduction}
the limit $L(f)$ is the same for both weights. Then
$W_1|f-L|=\frac{\varphi_1(h)}{\varphi_2(h)}\,W_2|f-L|$, so
$\sup_{\{h>R\}}W_1|f-L|\le M\sup_{t>R}\varphi_1(t)/\varphi_2(t)\to0$ uniformly in
$[f]\in\mathcal F$, which is hypothesis (iii) of \cref{thm:AA} for $\varphi_1$;
(i) follows from \cref{prop:global-comparison}, and (ii) because
$W_1(f-L)=(\varphi_1/\varphi_2)(h)\cdot W_2(f-L)$ is a product of a fixed continuous
function with an equicontinuous, uniformly bounded family. Apply \cref{thm:AA}.
\end{proof}

\begin{corollary}[Polynomial case]\label{cor:compact-poly}
If $p_2>p_1\ge0$ then $C^{\sharp}_{h;p_2}(A)\to C^{\sharp}_{h;p_1}(A)$ is compact on sets
that are bounded in $\|\cdot\|^{\sharp}_{\infty,h,p_2}$ and satisfy the equicontinuity
hypothesis of \cref{thm:compact-inclusion}.
\end{corollary}

\section{Measured extensions and \texorpdfstring{$L^q$}{Lq} embeddings}
\label{sec:Lq}

Let $\mu$ be a Radon measure on $A$ and $V(R):=\mu(\{h\le R\})$.

\begin{definition}[Polynomial volume growth]\label{def:volume}
$\mu$ has \emph{polynomial volume growth of order $\gamma\ge0$} if
$V(R)\le C_V(1+R)^{\gamma}$ for all $R\ge0$.
\end{definition}

\begin{theorem}[Sup-to-$L^q$ embedding, with loss]\label{thm:sup-to-Lq}
Assume $V(R)\le C_V(1+R)^{\gamma}$. Let $q\in[1,\infty)$, $p>p'\ge0$ and suppose
\[
(p-p')\,q\ >\ \gamma+1 .
\]
Then there is $C=C(p,p',q,\gamma,C_V)$ such that for all $f$ and $L$,
\[
\Bigl(\int_A\bigl((1+h)^{p'}|f-L|\bigr)^{q}d\mu\Bigr)^{1/q}
\ \le\ C\,\|f\|_{\infty,h,p;L}.
\]
\end{theorem}

\begin{proof}
Write $A_{-1}:=\{h\le0\}$ and $A_j:=\{j<h\le j+1\}$ for $j\in\N$, so that
$A=A_{-1}\sqcup\bigsqcup_{j\ge0}A_j$. On $A_{-1}$ the integrand is bounded by
$\|f\|_{\infty,h,p;L}^{q}$ and $\mu(A_{-1})\le V(0)\le C_V$, contributing at most
$C_V\|f\|_{\infty,h,p;L}^{q}$. On $A_j$, $j\ge0$, we have
$(1+h)^{p'}|f-L|\le\|f\|_{\infty,h,p;L}(1+h)^{p'-p}\le\|f\|_{\infty,h,p;L}(1+j)^{p'-p}$,
and $\mu(A_j)\le V(j+1)\le C_V(2+j)^{\gamma}\le C_V2^{\gamma}(1+j)^{\gamma}$. Hence
\[
\int_A\bigl((1+h)^{p'}|f-L|\bigr)^{q}d\mu
\le C_V2^{\gamma}\|f\|_{\infty,h,p;L}^{q}\sum_{j\ge0}(1+j)^{(p'-p)q+\gamma},
\]
and the series converges precisely when $(p-p')q-\gamma>1$.
\end{proof}

\begin{remark}[The loss is necessary]
The statement of the posted version had the \emph{same} weight $(1+h)^{p}$ on both sides
and the condition $pq>\gamma$. It is false: with $p'=p$ the layer estimate gives
$\sum_j\mu(A_j)=\mu(A)$, which is infinite as soon as $\mu$ is infinite. Concretely, on
$A=[0,\infty)$ with $\mu$ Lebesgue ($\gamma=1$), $h(x)=x$ and $f(x)=(1+x)^{-p}$, $L=0$,
one has $\|f\|_{\infty,h,p;0}=1$ while $\int_0^{\infty}((1+x)^{p}|f|)^{q}dx=\infty$.
\end{remark}

\begin{remark}
For a general $\varphi\in\Phi_{\mathrm{adm}}$ and a pair $\varphi'\preceq\varphi$, the same
layer decomposition gives the embedding under the summability hypothesis
$\sum_{j\ge0}(\varphi'(j)/\varphi(j))^{q}\bigl(V(j+1)-V(j)\bigr)<\infty$.
\end{remark}

\section{Functoriality under coarse-affine maps}
\label{sec:functoriality}

\begin{definition}[The category $\mathrm{Exh}$]\label{def:Exh}
Objects are pairs $(A,h)$ as in \cref{ass:standing}. A morphism
$\phi:(A,h_A)\to(B,h_B)$ is a continuous map $\phi:A\to B$ for which there are constants
$A_0\ge1$, $B_0\ge0$ with
\begin{equation}\label{eq:star}
h_A(x)\ \le\ A_0\,h_B(\phi(x))+B_0\qquad\text{for all }x\in A. \tag{$\star$}
\end{equation}
Composition and identities are the usual ones.
\end{definition}

\begin{remark}
Condition \eqref{eq:star} bounds $h_A$ \emph{above} in terms of $h_B\circ\phi$; it says
that $\phi$ does not send points that are far out in $A$ to points that are close in in
$B$. The contributions list of the posted version stated the inequality in the opposite
direction. Properness of $\phi$ is not needed for \cref{thm:pullback}.
\end{remark}

\begin{theorem}[Bounded pullback]\label{thm:pullback}
Let $\phi:(A,h_A)\to(B,h_B)$ satisfy \eqref{eq:star} and let $\varphi=\varphi_p$, $p>0$.
Then there is $C=C(A_0,B_0,p)$ such that for all $[f]\in C^{\sharp}_{h_B;\varphi}(B)$,
\[
\|\phi^{\ast}[f]\|^{\sharp}_{\infty,h_A;\varphi}\ \le\ C\,\|[f]\|^{\sharp}_{\infty,h_B;\varphi},
\qquad
C^{(h_A)}_{\varphi}(\phi^{\ast}[f])\ \le\ C\,C^{(h_B)}_{\varphi}([f]),
\]
where $\phi^{\ast}[f]:=[f\circ\phi]$. If in addition $h_B\circ\phi\le a_0h_A+b_0$ for some
$a_0\ge1$, $b_0\ge0$, and $\phi$ is surjective, then $\phi^{\ast}$ is a topological
isomorphism onto its image.
\end{theorem}

\begin{proof}
Let $L=L(f)$. By \eqref{eq:star},
$(1+h_A(x))^{p}\le(1+A_0h_B(\phi(x))+B_0)^{p}\le C_1(1+h_B(\phi(x)))^{p}$ with
$C_1=(\max(1,A_0)(1+B_0))^{p}$, using $1+A_0t+B_0\le\max(1,A_0)(1+B_0)(1+t)$. Hence
\[
\sup_{x\in A}(1+h_A(x))^{p}|f(\phi(x))-L|
\le C_1\sup_{x\in A}(1+h_B(\phi(x)))^{p}|f(\phi(x))-L|
\le C_1\|[f]\|^{\sharp}_{\infty,h_B;\varphi}.
\]
In particular $f\circ\phi$ is at bounded weighted distance from the constant $L$, so
$L(f\circ\phi)=L$ by \cref{thm:reduction} provided $h_A$ is unbounded on $A$, which holds
by properness; then the left-hand side is the sharp norm of $\phi^{\ast}[f]$. The
$\limsup$ version is identical. The last assertion follows by running the argument in the
reverse direction.
\end{proof}

\begin{proposition}[Contravariant functor]\label{prop:functor}
For fixed $\varphi=\varphi_p$, the assignment
$(A,h)\mapsto\bigl(C^{\sharp}_{h;\varphi}(A),\|\cdot\|^{\sharp}\bigr)$,
$\phi\mapsto\phi^{\ast}$, is a contravariant functor from $\mathrm{Exh}$ to Banach spaces
and bounded linear maps.
\end{proposition}

\section{Anisotropic spaces: ends and the correct quotient}
\label{sec:ends}

\subsection{The failure of the single-constant quotient}

\begin{assumption}[Decomposition at infinity]\label{ass:ends}
There is a compact $K_0\subset A$ (the \emph{core}) and finitely many pairwise disjoint,
open, connected, non-compact sets $U_1,\dots,U_m$ (the \emph{ends}) with
$A\setminus K_0=\bigsqcup_{i=1}^{m}U_i$. Each $U_i$ carries a continuous
$h_i:U_i\to[0,\infty)$ such that
\begin{enumerate}
\item[(E1)] $\{x\in U_i:h_i(x)\le R\}$ is relatively compact in $A$ for every $R<\infty$;
\item[(E2)] $h_i$ extends continuously to $\overline{U_i}$ with $h_i\equiv0$ on
$\partial U_i\subset K_0$;
\item[(E3)] $h_i$ is unbounded on $U_i$;
\item[(E4)] $\partial U_i=\overline{U_i}\setminus U_i\ne\varnothing$ (equivalently, no end
is clopen in $A$).
\end{enumerate}
Each end carries an admissible weight $\varphi_i\in\Phi_{\mathrm{adm}}$. We write
$\omega_i$ for the point at infinity of $U_i$: $a_k\to\omega_i$ means
$a_k\in U_i$ eventually and $h_i(a_k)\to\infty$.
\end{assumption}

Since the $U_i$ are connected, $h_i$ continuous with $h_i=0$ on $\partial U_i$ and
unbounded, $h_i$ takes every value in $[0,\infty)$.

Define the \emph{block weight} $W:A\to[1,\infty)$ by $W\equiv1$ on $K_0$ and
$W:=\varphi_i\circ h_i$ on $U_i$; it is continuous by (E2) and $\varphi_i(0)=1$.

\begin{proposition}[The naive quotient is useless]\label{prop:naive-fails}
Assume \cref{ass:ends} with all $\varphi_i$ unbounded, and let $f\in C(A)$ be such that
$\inf_{c\in\R}\sup_AW|f-c|<\infty$. Then $f$ has the \emph{same} limit at every end.
Consequently, if $f\to L_i$ at $\omega_i$ with $L_1\ne L_2$, then
$\inf_{c}\sup_AW|f-c|=+\infty$.
\end{proposition}

\begin{proof}
If $\sup_AW|f-c|=M<\infty$ then on each $U_i$, $|f-c|\le M/\varphi_i(h_i)\to0$ along
$\omega_i$, so $f\to c$ at every end.
\end{proof}

\begin{remark}
The posted version defined the anisotropic sharp norm as $\inf_{c\in\R}\sup_AW|f-c|$ and
then stated (Example 1 of its \S8.2) that for a strip with $L_-\ne L_+$ the asymptotic
constant is merely positive. By \cref{prop:naive-fails} both quantities are $+\infty$, so
the announced gluing lemma ``$\inf_c\sup_AW|f-c|\simeq\max\{\text{core},\max_i\text{end}\}$''
is false: the right-hand side is finite in that example and the left-hand side is not. The
proof given there also took an infimum over non-constant functions $c(x)$, which is not
the quantity being estimated.
\end{remark}

\begin{remark}[Why (E4) is needed]\label{rem:E4}
Without (E4) the functional of \cref{def:aniso-norm} is only a seminorm. If some $U_i$ is
clopen, then $f:=\mathbf{1}_{U_i}$ is continuous, constant on $K_0$ and on every $U_j$,
hence $N_0(f)=N_1(f)=\dots=N_m(f)=0$, while $[f]\ne0$ in $C(A)/\R$. Condition (E4) also
supplies the boundary point at which the end limits are tested in the proof of
\cref{thm:aniso-complete}. An earlier draft handled this in a parenthesis inside a proof;
it belongs in the hypotheses.
\end{remark}

\subsection{The end-by-end functional}

\begin{definition}[Anisotropic sharp norm]\label{def:aniso-norm}
Under \cref{ass:ends} put, for $f\in C(A)$,
\[
N_0(f):=\inf_{c\in\R}\ \sup_{x\in K_0}|f(x)-c|=\osc_{K_0}(f),
\qquad
N_i(f):=\inf_{c\in\R}\ \sup_{x\in U_i}\varphi_i(h_i(x))\,|f(x)-c|,
\]
and
\[
\|[f]\|^{\sharp}_{\mathbf{h},\boldsymbol{\varphi}}
:=\max\bigl\{N_0(f),\ N_1(f),\ \dots,\ N_m(f)\bigr\}\ \in[0,\infty].
\]
Let $\mathcal{A}:=\{[f]\in C(A)/\R:\ \|[f]\|^{\sharp}_{\mathbf{h},\boldsymbol{\varphi}}<\infty\}$.
\end{definition}

\begin{proposition}[It is a norm]\label{prop:aniso-is-norm}
$\|\cdot\|^{\sharp}_{\mathbf{h},\boldsymbol{\varphi}}$ is well defined on $C(A)/\R$
(invariant under adding constants), absolutely homogeneous, subadditive, and vanishes
only on the class of constants. Hence it is a norm on $\mathcal{A}$.
\end{proposition}

\begin{proof}
Each $N_j$ is invariant under $f\mapsto f+c$, homogeneous and subadditive (an infimum over
$c$ of a supremum of seminorms), and a maximum of seminorms is a seminorm. If
$\|[f]\|^{\sharp}=0$ then $f\equiv c_0$ on $K_0$ and $f\equiv c_i$ on $U_i$ for constants
$c_j$. By (E4) there is $y\in\partial U_i\subset K_0$; continuity of $f$ at $y$ forces
$c_i=c_0$. Hence $f$ is constant.
\end{proof}

\begin{theorem}[Anisotropic reduction]\label{thm:aniso-reduction}
Assume all $\varphi_i$ unbounded. Then $\|[f]\|^{\sharp}_{\mathbf{h},\boldsymbol{\varphi}}<\infty$
if and only if $f$ has a limit $L_i:=\lim_{a\to\omega_i}f(a)$ at each end and
$\varphi_i(h_i)\,|f-L_i|$ is bounded on $U_i$ for each $i$. In that case each infimum
$N_i$ is attained at the unique constant $c=L_i$, and
\[
\|[f]\|^{\sharp}_{\mathbf{h},\boldsymbol{\varphi}}
=\max\Bigl\{\osc_{K_0}(f),\ \max_{i}\ \sup_{U_i}\varphi_i(h_i)|f-L_i|\Bigr\}.
\]
\end{theorem}

\begin{proof}
Apply \cref{thm:reduction} on each end separately, with $A$ replaced by $U_i$ and $h$ by
$h_i$ (hypotheses (E1) and (E3) provide the properness and unboundedness used there);
$N_0(f)$ is finite automatically since $K_0$ is compact and $f$ continuous.
\end{proof}

\begin{theorem}[Completeness]\label{thm:aniso-complete}
$\bigl(\mathcal{A},\|\cdot\|^{\sharp}_{\mathbf{h},\boldsymbol{\varphi}}\bigr)$
is a Banach space.
\end{theorem}

\begin{proof}
Fix $x_0\in K_0$ and normalise representatives by $f_n(x_0)=0$. Let $([f_n])$ be Cauchy,
with $\|[f_n]\|^{\sharp}\le M$ for all $n$.

\emph{Uniform control.} $\osc_{K_0}(f_n)\le M$ and $f_n(x_0)=0$ give
$\sup_{K_0}|f_n|\le 2M$. Writing $L^{(n)}_i$ for the end-limits, continuity across
$\partial U_i$ and (E2) give, for $y\in\partial U_i$,
$|f_n(y)-L^{(n)}_i|\le\varphi_i(0)^{-1}\sup_{U_i}\varphi_i(h_i)|f_n-L^{(n)}_i|\le M$
(taking a limit of points of $U_i$ approaching $y$), hence $|L^{(n)}_i|\le3M$. Thus the
$f_n$ are uniformly bounded on $A$ and, on each $U_i$,
$g^{(n)}_i:=\varphi_i(h_i)(f_n-L^{(n)}_i)$ is uniformly bounded by $M$.

\emph{Passage to the limit.} The Cauchy property gives:
$\osc_{K_0}(f_n-f_k)\to0$, so $f_n-f_n(x_0)=f_n$ is uniformly Cauchy on $K_0$; and
$\sup_{U_i}\varphi_i(h_i)|(f_n-f_k)-(L^{(n)}_i-L^{(k)}_i)|\to0$, so
$(g^{(n)}_i)_n$ is uniformly Cauchy on $U_i$ and $(L^{(n)}_i)_n$ is Cauchy in $\R$ (test at
a point of $\partial U_i$ where $\varphi_i(h_i)=1$ and use the uniform convergence on
$K_0$). Let $f$ be the uniform limit on $K_0$, $L_i:=\lim_nL^{(n)}_i$,
$g_i:=\lim_ng^{(n)}_i$ uniformly on $U_i$, and set $f:=L_i+g_i/\varphi_i(h_i)$ on $U_i$.
The two definitions agree on $\partial U_i$ and $f$ is continuous. Finally
$\|[f_n]-[f]\|^{\sharp}=\max\{\osc_{K_0}(f_n-f),\max_i\sup_{U_i}|g^{(n)}_i-g_i|\}\to0$.
\end{proof}

\begin{proposition}[Vanishing-tail class]\label{prop:aniso-decay}
Let $\mathcal{A}_{\mathrm{van}}:=\{[f]\in\mathcal{A}:\ \limsup_{a\to\omega_i}
\varphi_i(h_i)|f-L_i|=0\ \text{for all }i\}$. Then $[f]\in\mathcal{A}_{\mathrm{van}}$ iff
$f\to L_i$ at each end with $f-L_i=o\bigl(\varphi_i(h_i)^{-1}\bigr)$, and
$\mathcal{A}_{\mathrm{van}}$ is a closed subspace of $\mathcal{A}$.
\end{proposition}

\begin{proof}
The characterisation is the definition. Closedness follows from
\cref{thm:aniso-complete} and the fact that a uniform limit of functions vanishing at
$\omega_i$ vanishes at $\omega_i$.
\end{proof}

\subsection{Two examples}

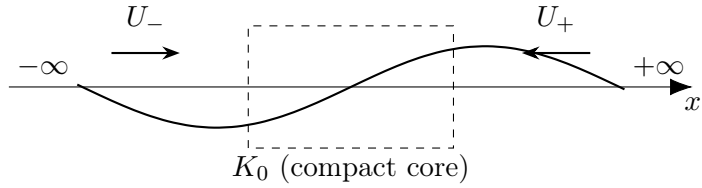
\begin{figure}[H]
\centering
\begin{tikzpicture}[scale=0.9,>=Stealth]
\draw[-{Latex[length=3mm]}] (-5,0)--(5,0) node[below]{$x$};
\draw (-4.5,0) node[above]{$-\infty$};
\draw (4.5,0) node[above]{$+\infty$};
\draw[thick,domain=-4:4,smooth] plot(\x,{0.6*sin(0.8*\x r)});
\node at (-3,1) {$U_{-}$};
\node at (3,1) {$U_{+}$};
\draw[dashed] (-1.5,-0.9) rectangle (1.5,0.9);
\node at (0,-1.2) {$K_0$ (compact core)};
\draw[->, thick] (-3.5,0.5) -- (-2.5,0.5);
\draw[->, thick] (3.5,0.5) -- (2.5,0.5);
\end{tikzpicture}
\caption{Two ends on a cylinder $\R\times M$: neighbourhoods $U_\pm$ and a compact core.}
\label{fig:cylinder}
\end{figure}

\begin{example}[Infinite strip]\label{ex:strip}
$A=\R\times[0,1]$, $K_0=[-1,1]\times[0,1]$, $U_\pm=\{\pm x>1\}$,
$h_\pm(x,y)=\max\{0,\pm x-1\}$, $\varphi_\pm(t)=(1+t)^{p}$. A continuous $f$ with
$f(x,y)\to L_\pm$ as $x\to\pm\infty$ at rate $O(|x|^{-p})$ has finite anisotropic norm
whether or not $L_-=L_+$. With the single-constant quotient of the posted version the norm
would be $+\infty$ as soon as $L_-\ne L_+$ (\cref{prop:naive-fails}); this is precisely the
defect the present formulation removes.
\end{example}

\begin{figure}[H]
\centering
\begin{tikzpicture}[scale=1.05,>=Stealth]
\draw[-{Latex[length=3mm]}] (-5,0)--(-0.1,0);
\draw[{Latex[length=3mm]}-] (0.1,0)--(5,0) node[below]{$x$};
\fill[white] (0,0) circle (0.12);
\draw[thick] (0,0) circle (0.12);
\node[below] at (0,-0.2) {$0$};
\draw (-4.5,0.3) node {$U_{-\infty}$};
\draw (4.5,0.3) node {$U_{+\infty}$};
\draw (-1.5,0.5) node {$U_{0^-}$};
\draw (1.5,0.5) node {$U_{0^+}$};
\draw[->, thick] (-4,0.2) -- (-3,0.2);
\draw[->, thick] (4,0.2) -- (3,0.2);
\draw[->, thick] (-0.8,0.5) -- (-0.2,0.5);
\draw[->, thick] (0.8,0.5) -- (0.2,0.5);
\draw[dashed] (-2.5,-0.3)--(-2.5,0.3);
\draw[dashed] (2.5,-0.3)--(2.5,0.3);
\node[below] at (-2.5,-0.4) {$-2$};
\node[below] at (2.5,-0.4) {$2$};
\end{tikzpicture}
\caption{Four ends on $\R\setminus\{0\}$: near $0^{\pm}$ and at $\pm\infty$.}
\label{fig:four-ends}
\end{figure}
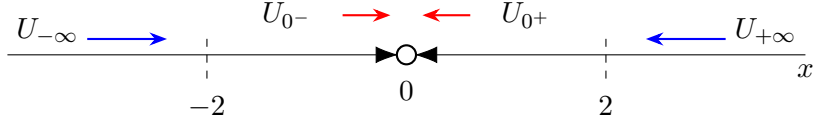

\begin{example}[Four ends on $\R\setminus\{0\}$]\label{ex:four-ends}
$A=\R\setminus\{0\}$, $K_0=[-2,-1]\cup[1,2]$,
$U_{\pm\infty}=\{\pm x>2\}$, $U_{0^{\pm}}=\{0<\pm x<1\}$, with
\[
h_{\pm\infty}(x)=\pm x-2,\qquad h_{0^{\pm}}(x)=-\log|x|,
\]
and $\varphi_{\pm\infty}(t)=(1+t)^{p}$, $\varphi_{0^{\pm}}(t)=e^{\alpha t}$. Note
$h_{0^{\pm}}$ satisfies (E2) since $h_{0^{\pm}}(\pm1)=0$, and
$\varphi_{0^{\pm}}(h_{0^{\pm}}(x))=|x|^{-\alpha}$. The anisotropic norm thus controls
polynomial decay at $\pm\infty$ and decay of order $|x|^{\alpha}$ towards $0^{\pm}$,
with four independent limits.
\end{example}

\subsection{What is left open in the anisotropic setting}

\begin{question}[Products]\label{q:products}
Given $(A,\{h^A_i\},\{\varphi^A_i\})$ and $(B,\{h^B_j\},\{\varphi^B_j\})$, is there a
canonical multi-exhaustion on $A\times B$ for which the anisotropic norm is submultiplicative
on simple tensors? The decomposition $A\times B=K^A_0\times K^B_0\ \sqcup\ \bigsqcup_{i,j}
U^A_i\times U^B_j$ proposed in the posted version is not a decomposition of $A\times B$ into
a compact core and open ends: the ``mixed'' pieces $K^A_0\times U^B_j$ and
$U^A_i\times K^B_0$ are omitted. A correct treatment must include them.
\end{question}

\begin{question}[Infinitely many ends]\label{q:infinite-ends}
Under a local finiteness assumption on the cover $\{U_i\}_{i\in I}$, does
\cref{thm:aniso-complete} persist with $\max_i$ replaced by $\sup_i$? Completeness of the
resulting space, and the correct topology on the family of limits $(L_i)_{i\in I}$
(bounded? converging to a common value along the end structure?), are not settled here.
\end{question}

\begin{question}[Complementation]\label{q:projection}
Is $\mathcal{A}_{\mathrm{van}}$ complemented in $\mathcal{A}$? With cut-offs $\chi_i$
supported near $\omega_i$ one may set $P[f]:=[f-\sum_iL_i\chi_i]$; this is linear
(the $L_i$ depend linearly on $f$) and idempotent, and $P[f]\in\mathcal{A}_{\mathrm{van}}$,
so the answer is affirmative once the boundedness $\|P[f]\|^{\sharp}\le C\|[f]\|^{\sharp}$
is established. That estimate requires control of $\osc_{K_0}(\sum_iL_i\chi_i)$ by
$\|[f]\|^{\sharp}$, hence a bound on $\max_i|L_i-L_j|$, which does \emph{not} follow from
the norm: the norm is invariant under adding a constant only, not under independent shifts
of the $L_i$. We leave the question open; the ``idea of construction'' given in the first
version did not address it.
\end{question}

\section{Conclusion and open problems}
\label{sec:conclusion}

Let us summarise what the exhaustion formalism does and does not deliver.

\paragraph{What it delivers.} A uniform notation for weighted decay across metric,
topological, ordered and measured settings; a clean statement of coarse-affine invariance
within a fixed class (\cref{lem:equiv-affine}, \cref{thm:coarse-stability}); functoriality
under coarse-affine maps (\cref{thm:pullback}); a corrected anisotropic formalism with a
reduction theorem and completeness (\cref{sec:ends}).

\paragraph{What it does not deliver.} A choice-free decay classification: by
\cref{ex:no-coarse-affine}, $O(h^{-p})$ depends on the coarse-affine class of $h$, and
distinct proper exhaustions of the same space can give genuinely distinct classes. Nor
does it deliver new function spaces in the isotropic case: by \cref{cor:isometry} the
sharp space \emph{is} $\Cb(A)$ and its decay class \emph{is} $\Cz(A)$, isometrically. Every
Banach-space statement about them is a classical statement in disguise, and should be
proved that way.

\paragraph{Where genuinely new questions remain.}

\begin{question}[Bounded weights]\label{q:bounded-phi}
For $\varphi$ bounded, $\|[f]\|^{\sharp}=\inf_c\sup_AW|f-c|$ is a weighted Chebyshev
approximation problem by constants on a non-compact space. Existence of minimisers (by
convexity and coercivity when $W$ is bounded below and $f$ is bounded), uniqueness,
alternation-type characterisations, and the structure of the contact set are all open in
this generality. Note that a genuine alternation criterion does hold: if the supremum
defining $\|[f]\|^{\sharp}$ is \emph{attained} at points $a_+,a_-$ with
$f(a_+)-c_\ast>0>f(a_-)-c_\ast$, then $c_\ast$ is the unique minimiser; the proof is that
moving $c$ in either direction strictly increases the contribution of one of the two
points. Attainment, however, is exactly what fails when $W$ is unbounded, and is not
automatic when $W$ is bounded on a non-compact space.
\end{question}

\begin{question}[Log-convexity in the bounded case]
See \cref{q:logconvex}.
\end{question}

\begin{question}[Genericity]
Is uniqueness of the minimising constant generic (in the Baire sense) in the bounded-weight
case, on the Banach space $\{f\in C(A):\ \|[f]\|^{\sharp}<\infty\}$ normed by
$|f(x_0)|+\|[f]\|^{\sharp}$? The argument sketched in the posted version establishes neither
that the exceptional sets are closed nor that they have empty interior.
\end{question}

\begin{question}[Anisotropic duality and compactness]
Identify $\mathcal{A}^{\ast}$ and $\mathcal{A}_{\mathrm{van}}^{\ast}$, and formulate the
analogue of \cref{thm:AA}. The natural guess, by \cref{thm:aniso-reduction}, is a closed
subspace of $C(K_0)/\R\oplus_{\infty}\bigoplus_{i,\infty}\Cb(U_i)$ cut out by the matching
conditions on $\partial U_i$; making this precise, and describing the dual, is not done
here.
\end{question}

\begin{question}[An application that uses the anisotropy]
The formalism of \cref{sec:ends} is designed for elliptic problems on manifolds with ends,
where solutions decay at different rates in different channels. Producing a statement that
is not already contained in the $b$-calculus \cite{Melrose1993} or in
Kondrat'ev's theory \cite{Kondratev1967} is the real test, and remains to be done.
\end{question}

\appendix

\section{Corrections}
\label{app:corrections}

This appendix has two parts. Part~\ref{app:part-posted} lists the statements of the
\emph{posted} version of this work that are retracted or substantially modified here;
these are the corrections a reader comparing the two documents will want.
Part~\ref{app:part-revision} lists errors that were \emph{not} present in the posted
version: they were introduced in unposted drafts during this revision and corrected before
posting. The second list is included because the mistakes are instructive, because the
same failure mode occurs twice, and because it would be dishonest to present a corrected
paper as having been arrived at without incident.

\subsection{Statements retracted from the posted version}
\label{app:part-posted}

The following appeared in the posted version and are retracted or substantially modified.
We list them in the order in which they occurred there.

\begin{enumerate}

\item \textbf{Coarse affine equivalence of proper exhaustions (Prop.~2.3).}
\emph{False.} Counterexample: $A=\R$, $h(x)=|x|$, $h'(x)=x^{2}$
(\cref{ex:no-coarse-affine}). The proof deduced affine growth from the mutual
quasi-inverse relations $R\le f(g(R))$, $R'\le g(f(R'))$, which do not imply it. As a
consequence, the accompanying Remark (``regularity of $h$ does not matter, all proper
exhaustions being coarsely equivalent'') is withdrawn, and coarse-affine equivalence is
now a standing hypothesis (\cref{ass:standing}).

\item \textbf{Universal property (Thm.~3.5).} \emph{Quantifiers inverted.} As stated the
hypothesis is satisfied by every $F$ (take $R=0$, $U=Y$) and the conclusion fails.
Corrected as \cref{thm:universal}.

\item \textbf{Definition of the spaces.} $C^{\sharp}_{h,p}(A):=C(A)/\R$ is not a normed
space: the sharp norm is $+\infty$ on most classes. All statements are now made on the
finiteness domain (\cref{rem:space-def}).

\item \textbf{Orlicz sharp norms (Def.~4.4, Thm.~4.5).} \emph{Removed.} The proof of the
norm axioms and completeness consisted of unjustified assertions (``$\le$ convex
combination $\le1$'', ``the $\Delta_2$ condition implies uniform convergence on
compacts''). The construction may well be sound but is not established, and by
\cref{thm:reduction} it would in any case reduce to an Orlicz-type renorming of $\Cb(A)$.

\item \textbf{Patching principle and exact patch formula (Lem.~4.7, Prop.~4.9).}
\emph{False as stated.} The local and tail functionals take independent infima over the
constant, so their maximum need not dominate the sharp norm. Replaced by
\cref{lem:splitting} (splitting at a fixed constant), which is what all downstream
arguments actually use.

\item \textbf{Tail control implies finite global norm (Thm.~4.11).} \emph{Vacuous.}
$R\mapsto Q_p(R)$ is non-increasing by construction, so ``eventually non-increasing'' is
automatic and the conclusion is a restatement of the definition. Superseded by
\cref{cor:finiteness}.

\item \textbf{Banach lattice structure (Prop.~4.12).} \emph{False and ill-posed.} The
order on $C(A)/\R$ is not well defined, and
$\|\,|[f]|\,\|^{\sharp}=\|[f]\|^{\sharp}$ fails (\cref{rem:no-lattice}). Only the
contraction inequality survives.

\item \textbf{Existence of a sharp minimiser (Prop.~4.13).} \emph{Correct conclusion,
defective proof, and superseded.} The lower semicontinuity argument passes to the limit
$R\to\infty$ in an inequality whose right-hand side depends on $R$ through
$\max_{h\le R}\varphi(h)\to\infty$. For unbounded $\varphi$ the statement is trivial by
\cref{thm:reduction}; for bounded $\varphi$ it is part of \cref{q:bounded-phi}.

\item \textbf{Two-sided contact condition (Lem.~4.14).} \emph{False.} Take $f\ge L$ with
$f\not\equiv L$: the contact set has no negative branch. The stated proof is absent.

\item \textbf{Baire-generic uniqueness (Thm.~4.15).} \emph{Retracted.} For unbounded
$\varphi$ the minimiser is always unique (\cref{thm:reduction}), so the statement is
vacuous; for bounded $\varphi$ neither closedness nor emptiness of interior of the sets
$E_k$ is established (``a small perturbation would break the tie'' is not an argument).
See \cref{q:bounded-phi}.

\item \textbf{Stability of minimisers (Prop.~4.16).} \emph{Ill-posed.} $J_f$ depends on the
representative $f$, not on the class $[f]$, so convergence in the quotient norm does not
control the minimisers. A correct version requires $\sup_AW|f_n-f|\to0$.

\item \textbf{Uniqueness criterion (Prop.~4.17).} \emph{Correct, with a simpler proof and
without the plateau hypothesis}, provided the supremum is \emph{attained}; see
\cref{q:bounded-phi}.

\item \textbf{Real interpolation (Thm.~4.19).} \emph{Retracted.} Its proof invokes the
log-convexity of Thm.~4.24, which is not established (item 15), and the ``nearly optimal
decomposition by truncating at level $R(t)$'' is not constructed.

\item \textbf{Global embeddings modulo a local patch (Prop.~4.21).} \emph{Degenerate.}
The quantity $\osc_{[f]}(a)=\inf_c|f(a)-c|$ is identically $0$. Replaced by
\cref{prop:global-comparison}.

\item \textbf{Log-convexity in $p$ (Thm.~4.24).} \emph{Not proved as stated.} H\"older
gives log-convexity of $p\mapsto\|f\|_{\infty,h,p;c}$ at fixed $c$; the infimum of a family
of log-convex functions need not be log-convex, and the three infima are attained at
different constants. For unbounded $\varphi$ the statement is nevertheless true, because
all three infima are attained at the same $c=L$ (\cref{prop:logconvex}); for bounded
$\varphi$ it is open.

\item \textbf{Moreau regularisation (Prop.~4.25).} \emph{Retracted.} The proof bounds the
inf-convolution by an expression containing $\exp(-d(x,y)^{2}/2\lambda)$, which does not
occur in the definition of $\mathsf{M}_\lambda f$; the argument is one for a Gaussian
kernel, not for a Moreau envelope. Moreover $\mathsf{M}_\lambda f\equiv-\infty$ unless $f$
is bounded below, and for $f$ Lipschitz the classical bound
$0\le f-\mathsf{M}_\lambda f\le\lambda\Lip(f)^{2}/2$ is uniform and therefore does
\emph{not} give convergence in a norm with unbounded weight.

\item \textbf{Weighted Schur test (Thm.~4.36).} \emph{Wrong normalisation.} The second
Schur condition is never used, and $\int K(x,y)d\mu(x)=1$ should be
$\int K(x,y)d\mu(y)=1$, which is what makes $T_K$ preserve constants and descend to the
quotient. Corrected as \cref{thm:schur}.

\item \textbf{Isometry with $\Cz$ (Thm.~4.37).} \emph{Correct conclusion, defective proof,
and strengthened.} Linearity was asserted with the admission that ``the situation is more
subtle''. For unbounded $\varphi$ the selection $c_\ast([f])=L(f)$ is canonical and linear,
and the map is a surjective isometry onto $\Cb(A)$, restricting to $\Cz(A)$ on the decay
class (\cref{cor:isometry}).

\item \textbf{Full duality (Thm.~5.3) and partial representation (Thm.~5.2).}
\emph{False.} $\Msp_0^{W^{-1}}(A)$ embeds isometrically with dense range but is not the
whole dual; $\Xdec^{\ast}\cong\Msp(A)$ (\cref{thm:duality}). The proofs were also
circular: Thm.~5.3 invoked Prop.~5.7 and Prop.~5.7 invoked Thm.~5.3.

\item \textbf{Optimality certificate (Cor.~5.4) and existence of dual extremals
(Cor.~5.8).} \emph{False.} \Cref{ex:no-extremal} exhibits $[f]$ whose sharp norm is not
attained on the unit ball of $\Msp_0^{W^{-1}}(A)$.

\item \textbf{Subdifferential description (Thm.~5.6).} \emph{False.} The normalisation
$\nu(\mathcal{C}^{+}_\ast)=-\nu(\mathcal{C}^{-}_\ast)=\tfrac12$ is not derived and confuses
raw with weighted mass; the contact set need not meet both signs; and the result depends on
the retracted Thm.~5.3.

\item \textbf{Weak-$*$ compactness (Prop.~5.7).} \emph{False}, being an application of
Banach--Alaoglu to a space that is not the dual in question.

\item \textbf{Weighted Arzel\`a--Ascoli (Thm.~5.10).} \emph{Correct in substance, with a
misplaced hypothesis.} Equicontinuity must be assumed for $W(f-L(f))$, not for $f$.
Corrected and reproved in three lines via \cref{cor:isometry} as \cref{thm:AA}.

\item \textbf{Sup-to-$L^q$ embedding (Thm.~6.3).} \emph{False.} The same weight cannot
appear on both sides; a strict loss $p>p'$ with $(p-p')q>\gamma+1$ is needed. The proof
also used $\mu(A_R)\lesssim(1+R)^{\gamma-1}$, which does not follow from a bound on $V$.
Corrected as \cref{thm:sup-to-Lq}.

\item \textbf{Direction of the functoriality condition.} The contributions list stated
$h_A\ge c\,h_B\circ\phi-C$; the condition used, and the one that makes the pullback
bounded, is $h_A\le A_0\,h_B\circ\phi+B_0$.

\item \textbf{Anisotropic sharp norm and gluing (Def.~8.5, Lem.~8.6).} \emph{False.}
The single-constant infimum is $+\infty$ whenever two ends carry different limits
(\cref{prop:naive-fails}), contradicting the announced equivalence; and Step 2 of the proof
took an infimum over non-constant functions $c(x)$, which is not the quantity defined.
Replaced by \cref{def:aniso-norm} and \cref{thm:aniso-reduction}.

\item \textbf{Products (Prop.~8.9), infinitely many ends (Prop.~8.10), bounded projection
(Thm.~8.11).} \emph{Retracted}, and restated as \cref{q:products},
\cref{q:infinite-ends}, \cref{q:projection}. The proposed product decomposition omits the
mixed pieces; the projection estimate requires a control on $\max_i|L_i-L_j|$ that the norm
does not provide.

\item \textbf{Conclusion.} The posted version described the framework as ``intrinsically
isotropic'' and proposed multi-exhaustions as future work, while \S8 of the same document
carried them out. The inconsistency is removed.

\item \textbf{Companion note on Berry--Esseen bounds.} The claim, made in the
introduction and in the ``Potential Applications'' section, that the weighted metric
recovers the rate $n^{-1/2}$ under $\mathbb{E}|X|^{2+\delta}<\infty$ is withdrawn; see
\cref{sec:intro}.

\end{enumerate}

\subsection{Errors made and corrected during the revision}
\label{app:part-revision}

None of the following was present in the posted version; all were introduced in unposted
drafts and are corrected in the present text.

\begin{enumerate}

\item \textbf{Coarse stability, first attempt.} \emph{False.} Counterexample:
$A=[0,\infty)$, $h(x)=x$, $h'(x)=2x$, $\varphi(t)=e^{t}$, $f-L=e^{-x}$, for which
$C^{(h)}_\varphi=1$ and $C^{(h')}_\varphi=+\infty$ (\cref{ex:exp-not-invariant}). The
hypothesis ``power-dilation control'' with exponent $\kappa>1$ is too weak, and the
exponential weights, explicitly claimed there to satisfy the conclusion, are exactly the
ones for which it fails. The invalid step is
$W^{\kappa}|f-L|\le(W|f-L|)^{\kappa}$, which is equivalent to $|f-L|\ge1$; see
\cref{rem:invariance-history}. Replaced by the equivalence \cref{thm:coarse-stability}, with
\cref{cor:poly-invariance} covering the polynomial and log-polynomial cases.

\item \textbf{Cascading wording.} The standing \cref{ass:standing} and point (1) of the
introduction asserted invariance of the norms and decay classes under coarse-affine
changes of $h$ without restriction on the weight. Both are now restricted to weights of
bounded dilation. \Cref{lem:equiv-affine} and \cref{cor:coarse-iso}, being stated only for
$\varphi_p$, are unaffected.

\item \textbf{Example of a non-attained supremum.} \emph{Stated in the wrong space.} The function $f(x)=(1+x)^{-1}$ gives $T([f])\equiv1\in\Cb\setminus\Cz$, so
$[f]\notin\Xdec$ and the example did not instantiate the corollary it was attached to.
Replaced by $f(x)=(1+x)^{-2}$ in \cref{ex:no-extremal}; the conclusion is unchanged.

\item \textbf{Existence of extremals.} Non-attainment is not systematic. If the contact
set of $T([f])$ meets both signs, a norming measure exists
(\cref{prop:extremal-iff}); the posted version did not distinguish the two cases.

\item \textbf{Weighted Schur test.} The continuity of $T_Kf$ was taken
for granted; a hypothesis to that effect is now stated in \cref{thm:schur}.

\item \textbf{Anisotropic hypotheses.} The
non-emptiness of $\partial U_i$ was handled by a parenthetical remark inside a proof. It is
now hypothesis (E4) of \cref{ass:ends}; without it $f=\mathbf{1}_{U_i}$ has zero norm
without being constant (\cref{rem:E4}).

\item \textbf{Equivalence of exhaustion functions.} The proof used
$a\in K_{h(a)}$, i.e.\ attainment of the infimum in \cref{def:exhaustion} --- the very
caveat flagged elsewhere in the same section. Corrected with an $\varepsilon$-argument.

\item \textbf{Interpolation for asymptotic constants.} The
indeterminate case $0\cdot\infty$ was not excluded; $u_n=1/n$, $v_n=n$ shows the
inequality then fails. A hypothesis has been added.

\item \textbf{Sup-to-$L^q$ embedding.} The layer $A_{-1}=\{h\le0\}$ was
defined but omitted from the summation. Its (trivial) contribution is now included.

\item \textbf{Log-convexity.} Posed as an open question and
then resolved three lines later in the unbounded case. Split into
\cref{prop:logconvex} (proved, unbounded weights) and \cref{q:logconvex} (open, bounded
weights).

\item \textbf{Coarse-affine invariance, second attempt.} \emph{False as an equivalence.} The proof of (iii)$\Rightarrow$(i) selected $a_n\in A$ with
$h(a_n)=t_n$, which presupposes that $h$ attains all large values --- true when $A$ is
connected, false in general, and false for the ``Directed'' instance $A=\N$ listed in
\cref{ex:instances} of the paper itself. Worse, the statement fails:
\cref{ex:N-inv-strict} gives a weight on $A=\N$ with unbounded global dilation for which
the asymptotic constants are nevertheless invariant under every coarse-affine change of
$h$.

\item \textbf{Relativising to $h(A)$ is not sufficient either.} The natural repair ---
replace $\sup_{t\ge0}$ by $\sup_{t\in h(A)}$ --- restores the necessity argument but breaks
sufficiency: the proof of the lower bound $W\le CW'$ requires dilation control at the
points of $h'(A)$, not of $h(A)$, and these differ. \Cref{ex:N-van-strict} exhibits
$\varphi$ on $A=\N$ with $D_h$ finite for which comparability fails, via $h'=h/2$.
\Cref{thm:coarse-stability} therefore separates the two invariance properties and
characterises each by its own dilation functional: comparability by $D_h$ \emph{and}
$E_h$, preservation of vanishing by $D_h$ alone. Part (c) recovers the discarded statement
under the hypothesis $h(A)\supseteq[t_0,\infty)$, which holds for connected $A$.

\item \textbf{Existence of norming measures.} An intermediate draft left necessity open. It
is necessary, and the criterion is exact: a norming measure exists if and only if
$g=T([f])$ attains both $+\|g\|_\infty$ and $-\|g\|_\infty$
(\cref{prop:extremal-iff}). The proof is a Jordan decomposition plus the observation that
a one-signed $\sigma$ cannot satisfy $\int W\,d\sigma=0$ when $W\ge1$.

\item \textbf{Weighted Schur test.} The parenthetical justification of the continuity of
$T_Kf$ referred to differentiation under the integral sign; it should refer to passing to
the limit under the integral sign (dominated convergence).

\item \textbf{Choice of the sequence $t_n$.} In the necessity argument the $t_n$ must be
strictly increasing for the interpolation defining $v$ to be legitimate; a subsequence is
now extracted explicitly.

\end{enumerate}

\end{document}